\documentclass[12pt]{amsart}

\usepackage{epsfig}

\newcommand{\tr}{{\rm tr}}
\newcommand{\cL}{{\mathcal L}}
\newcommand{\cM}{{\mathcal M}}
\newcommand{\cG}{{\mathcal G}}
\newcommand{\cA}{{\mathcal A}}
\newcommand{\cW}{{\mathcal W}}
\newcommand{\cX}{{\mathcal X}}
\newcommand{\bA}{{\hspace{-0.3pt}\mathbb A}\hspace{0.3pt}}
\newcommand{\bC}{{\hspace{-0.3pt}\mathbb C}\hspace{0.3pt}}
\newcommand{\Hom}{\operatorname{Hom}}

\newcommand{\smlst}[1]{\makebox[0pt][r]{\scriptsize{$#1$}}}
\newcommand{\smrst}[1]{\makebox[0pt][l]{\scriptsize{$#1$}}}
\newcommand{\smst}[1]{\makebox[0pt]{\scriptsize{$#1$}}}

\newcommand{\g}{SL_2(\bC)} 
\newcommand{\lmk}{\noindent\mbox{}\hfill}
\newcommand{\rmk}{\hfill\mbox{}\par}

\newcommand{\lcr}{\raisebox{-5pt}{\mbox{}\hspace{1pt}
                  \epsfig{file=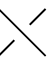}\hspace{1pt}\mbox{}}}
\newcommand{\ift}{\raisebox{-5pt}{\mbox{}\hspace{1pt}
                  \epsfig{file=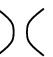}\hspace{1pt}\mbox{}}}
\newcommand{\zer}{\raisebox{-5pt}{\mbox{}\hspace{1pt}
                  \epsfig{file=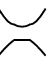}\hspace{1pt}\mbox{}}}

\newtheorem{theorem}{Theorem}
\newtheorem{proposition}{Proposition}
\newtheorem{definition}{Definition}
\newtheorem{lemma}{Lemma}
\newtheorem{corollary}{Corollary}

\newcommand{\U}{U_h(sl_2)}
\newcommand{\V}{\underline{\frac{1}{2}}}

\title{Topological Interpretations of Lattice Gauge Field Theory}
\author{doug bullock} 
\address{Department of Mathematics, Boise State University, Boise, ID
83725, USA}
\email{bullock@math.idbsu.edu}
\author{CHARLES FROHMAN}
\address{Department of Mathematics, University of Iowa, Iowa City, IA
52242, USA}
\email{frohman@math.uiowa.edu}
\author{JOANNA  KANIA-BARTOSZY\'{N}SKA}
\address{Department of Mathematics, Boise State University, Boise, ID
83725, USA}
\email{kania@math.idbsu.edu}
\thanks{The first author is partially supported by an Idaho SBOE
Specific Research Grant; the second and third by NSF-DMS-9204489 and
NSF-DMS-9626818.}

\begin{document}

\begin{abstract} We construct lattice gauge field theory based on a
quantum group 
on a lattice of dimension 1. Innovations include a coalgebra structure
on the connections, and an investigation of connections that are not
distinguishable by observables. We prove that when the quantum group
is a deformation of a connected algebraic group (over the complex numbers), 
then the algebra of observables forms a deformation quantization of the ring 
of characters of the fundamental group of the lattice with respect to the 
corresponding algebraic group.  Finally, we investigate lattice gauge field 
theory based on quantum $SL_2{\mathbb C}$, and conclude that the algebra of 
observables 
is the Kauffman bracket skein module of a cylinder over a surface associated 
to the lattice.\end{abstract}

\maketitle

\part*{Introduction}

Lattice gauge field theory based on an algebraic group $G$ is a finite
element approximation of a smooth gauge field theory with $G$ as its
structure group.  Infinitesimally varying connections and gauge
transformations on a principal bundle are discretized via a lattice
embedded in the base manifold.  To each edge in the lattice a
connection imparts an element of $G$ encoding the holonomy along that
edge.  Gauge fields (i.e., functions on connections) are represented
by a copy of the coordinate ring of $G$ associated to each edge.  The
action of the gauge group is then concentrated at vertices.  All
computations become merely algebraic with analytic and geometric
considerations swept aside. The end result is an algebra of
observables (the gauge invariant gauge fields), that can be
understood as the character theory for representations of the
fundamental group of the lattice into $G$.

Lattice gauge field theory based on a quantum group yields a
deformation of this theory.  Technically, the result is an algebra of
observables that, with respect to the standard Poisson structure
\cite{BG,Go}, gives a deformation quantization of the ring of
$G$-characters.  Here one must think of the fundamental group of the
lattice as the fundamental group of a surface with boundary.  

In this paper we develop, from an elementary and computational
viewpoint, the basic objects of lattice gauge field theory based on a
ribbon Hopf algebra, of which a quantum group is an example.  We then
use this foundation to begin the study of the structure of algebras of
observables (paying particular attention to quantum groups) and to
recognize the observables as algebras that have already been studied
in a topological framework.

The genesis of our approach can be found in the papers
\cite{AGS,Bu,BR1,BR2,FR}.  Fock and Rosly were the first to derive the
Poisson structure on $G$-characters from a lattice gauge field theory.
Their formula is written in terms of a solution of the modified
classical Yang-Baxter equation \cite{drinfeld}.  Also, recall that the
characters of a surface group are only a homotopy invariant while the
Poisson structure is a topological invariant.  For this reason, Fock
and Rosly endow the lattice with extra information, called a
ciliation, so that it determines a surface.

Passing to quantum groups, Alekseev, Grosse and Schomerus defined an
exchange algebra over a ciliated lattice so that basic elements of the
algebra of gauge fields commute according to a solution of the quantum
Yang-Baxter equation.  This algebra is related to a quantization of
the characters with respect to the usual Poisson structure.  It is
based on a full solution of the Clebsch-Gordan problem for the quantum
group being used, and has both gauge transformations and gauge fields
in the same place.

Buffenoir and Roche \cite{BR1} took this approach farther. First they
isolated the gauge fields from the gauge transformations.  Their gauge
algebra is dual to that of \cite{AGS}, hence they have a coaction of
the gauge algebra on the gauge fields.  The coinvariant part of the
gauge fields is the algebra of observables, which is a deformation of
the classical ring of characters.  They proceed to define Wilson loops
and the Yang-Mills measure and to derive $3$-manifold invariants from
this setting \cite{Bu,BR2}.

We found ourselves unable to compute examples in the exchange algebra
formulation. We instead define our gauge fields as ``functions'' on
the space of connections.  This  makes the structure of the algebra of
observables more clear.  Working from the point of view of
low-dimensional topology, we assume a familiarity with the
basics of knot theory.  Otherwise, one can read most of
this paper knowing only the definition of a ribbon Hopf algebra and a
smattering of its representation theory.  Kassel \cite{kassel} and
Sweedler \cite{sweedler} are sufficient references.

Part 1 is devoted to our translation of the
basic objects of a lattice gauge field theory and to our devices for
computing in the reformulated version.  We do not merely alter the
language of \cite{BR1}; there are three significant innovations which
provide the added computing power. The first is to realize that gauge
fields come from the restricted dual of the Hopf algebra on which the
theory is based.  This leads to a coordinate free formulation. Next,
we do not multiply gauge fields as abstract variables modulo exchange
relations.  Rather we comultiply connections in a way that implies the
usual exchange relations for fields while preserving their
evaluabililty.  Finally, we are able to mimic the classical phenomenon
of pushing the support of a gauge field around.  Our new foundations
allow us to compute Wilson loops and many other operators using a
simple extension of tangle functors.

The second part is devoted to an analysis of the structure of the
algebra of observables.  Our viewpoint is that the observables
corresponding to quantum groups generalize the rings studied by Procesi
\cite{procesi}.  He arrived at these rings as the invariants of
n-tuples of matrices under conjugation.  The connection with lattice
gauge field theory is that each n-tuple of matrices corresponds to a
connection on a lattice with one vertex and n-edges, with the gauge
fields  based on a classical group.

In passing from Procesi's work to ours, we find that the algebra of
observables corresponding to a quantum group is a more subtle object.
Instead of depending solely on the fundamental group of the lattice, the
observables are classified by the topological type of a surface
specified by a ciliated lattice. The construction given in this paper
leads to an algebra of ``characters'' of a surface group with respect
to any ribbon Hopf algebra.  The algebras are interesting from
many points of view: They generalize objects studied in invariant
theory; they should provide tools for investigating the structure of the
mapping class groups of surfaces; and they should give a way of
understanding quantum invariants of $3$-manifolds.

In the case that the data correspond to a connected affine algebraic
group $G$, it is possible to make explicit parallels with the existing
theory. The algebra of observables based on $U(\mathfrak{g})$ is
proved to be the ring of $G$-characters of the fundamental group of
the associated surface. Then, the original motivating problem is
solved: Given the ring of $G$-characters of a surface group, show that
the observables based on the corresponding Drinfeld-Jimbo algebra
form a quantization with respect to the usual Poisson structure.  We
also prove for the classical groups that the algebra of observables is
generated by Wilson loops.  Finally, invoking a quantized
Cayley-Hamilton identity, we obtain a new proof, independent of
\cite{quant}, that the $U_h(sl_2)$-characters of a surface are exactly
the Kauffman bracket skein module of a cylinder over that surface.

Many further avenues of research present themselves.  Working with
quantum groups defined over local fields side steps several
interesting and subtle structural questions. What happens when one uses
a quantum group at a root of unity?  How about lattice gauge field
theory based on a quasitriangular quasi-Hopf algebra?

There is a graphical calculus of the characters of the fundamental
group of any manifold with respect to any algebraic group, for
example, see \cite{isomorphism}, \cite{PS2} or \cite{FS}.  It derives
from the fact that Wilson loops are a pictorial description of
characters of the fundamental group of a manifold, after which the
tools of classical invariant theory express all functional relations
between characters in a diagrammatic fashion.  The graphical models
have only been worked out for special linear groups.

The power of lattice gauge field theory is that it places the
representation theory of the underlying manifold and the quantum
invariants in the same setting.  Ultimately the asymptotic analysis of
the quantum invariants of a $3$-manifold in terms of the
representations of its fundamental group should flow out of this
setting. The identification of the representation theory of a quantum
group with that of a compact Lie group leads to rigorous integral
formulas for quantum invariants of $3$-manifolds. This should in turn
lead to a simple explication of the relationship between quantum
invariants and more classical invariants of $3$-manifolds.

Finally, there should be a similarly clean treatment of lattice gauge
field theory for lattices of higher dimension. Although, in dimension
greater than two, the answers will no longer be topological in nature,
the constructions and objects should be of interest to
geometers, algebraists and analysts.

This work was done at the Banach Center, the University of Iowa, Boise
State University, The George Washington University, the University of
Missouri, and MSRI. The authors thank all their hosts for their
hospitality. We also thank Daniel Altschuler, Jorgen Anderson, Georgia
Benkhart, Vic Camillo, Fred Goodman, Joseph Mattes, Michael Polyak,
Florin Radulescu, Arun Ram, Justin Roberts, Don Schack and Bob Sulanke
for helpful conversations.

\part{Lattice Gauge Field Theory}

Herein we develop, from a self contained and axiomatic approach, the
machinery of gauge field theory on an abstract, oriented, ciliated
graph.  For basic background on Hopf algebras we rely on Sweedler
\cite{sweedler} and Kassel \cite{kassel}.  The discussion here is
restricted to the definitions and basic results confirming that the
theory is consistent and computationally viable.  For origins of the
ideas we refer the reader to \cite{AGS,Bu,BR1,BR2,survey,FR}.  

\section{Objects}

The elementary objects of a lattice gauge field theory are: a ribbon
Hopf algebra; an abstract graph which is oriented and ciliated; 
discretized connections, gauge transformations, and gauge fields.

\subsection{}\label{hopfalgebras}

Let $H$ be a ribbon Hopf algebra defined over a field $\mathfrak{k}$
or its power series ring $\mathfrak{k}[[h]]$.  In the latter case all
objects carry the h-adic topology (see \cite{kassel}), and all
morphisms are continuous.  In discussions germane to both settings we
will refer to the base over which the algebra is defined as
$\mathfrak{b}$.  Following Kassel we let $\mu$, $\eta$, $\Delta$,
$\epsilon$ and $S$ denote the multiplication, unit, comultiplication,
counit and antipode of $H$.  The universal $R$-matrix is $R = \sum_i
s_i \otimes t_i$, which we usually write as $s \otimes t$ with
summation understood.  The ribbon element is $\theta$, and we add a
{\bf charmed element}, $k = \theta^{-1} S(t)s.$ The charmed element is
{\bf grouplike}, meaning $\Delta(k) = k \otimes k$, and it satisfies
$k^{-1} = S(k) = \theta t S^2(s)$ and $S^2(x) = kxk^{-1}$ for all $x
\in H$.

The Hopf algebra dual is well documented in \cite{sweedler} provided
$\mathfrak{b} = \mathfrak{k}$.  The topological case, however, has
been neglected.  If $H$ is a Hopf algebra over $\mathfrak{k}[[h]]$
then the sets $U_n = \{L \in H^*\;|\; L(H) \subset h^n
\mathfrak{k}[[h]] \}$ form a neighborhood basis of the origin.  An
ideal $J \leq H$ is cofinite if $H/J$ is topologically free and
modeled on a finite dimensional vector space; $L \in H^*$ is cofinite
if $\ker(L)$ contains a cofinite ideal.  The restricted dual, $H^o$,
is the completion of the cofinite functionals.

It is not hard to check that $H^o$ is topologically free. In the case
that $H$ is a Drinfeld-Jimbo deformation of a simple Lie algebra
\cite{kassel}, it is clear that $H^o$ is modeled on $(H/hH)^o$.  To
see that it is a Hopf algebra one must check that $\mu^*$, $\eta^*$,
$\Delta^*$, $\epsilon^*$ and $S^*$ restricted to $H^o$ or $H^o \otimes
H^o$ take values in the appropriate spaces.  The only point that needs
any discussion is why $\mu^*(L) \in H^o \otimes H^o$.  Suppose that
$L$ is the limit of cofinite linear functionals $\{L_i\}$.  It follows
that $\{L_i\}$ is a Cauchy sequence in $H^o$.  The classical
discussion of $\mu^*$ in \cite{sweedler} provides $L_i \circ \mu \in
H^o \otimes H^o$, which is complete by definition \cite{kassel}. It is
easy to see that $L_i \circ \mu$ is also Cauchy, so continuity of
$\mu^*$ gives $\mu^*(L) \in H^o \otimes H^o$.

The Hopf algebra $H$ acts on its restricted dual in two obvious ways:
$x \cdot \phi(y)=\phi(xy)$ and $x \cdot \phi(y)=\phi(yx)$.  A
subalgebra of $H^{\circ}$ is stable if it is invariant under both
actions.  For the remainder of this section we fix a stable subalgebra
$B$ of $H^{\circ}$.  (Stability implies that $B$ is actually a Hopf
subalgebra.)

The adjoint action $ad : H \otimes H \rightarrow H$ is given, in
Sweedler notation, by \[ ad(Z)W = \sum_{(Z)} Z''XS(Z'),\] which is
further compressed to $ad(Z)W = Z''XS(Z')$.  (Readers unfamiliar with
this notation for comultiplication are refered to \cite{kassel}.)  By
taking duals we get the adjoint action of $H$ on $B$: if $Z,X \in H$,
and $\phi \in B$ then \[ ad(Z)(\phi)X= \phi(Z''XS(Z')).\] An element
$\phi \in B$ is {\bf invariant} if for every $Z \in H$, $ad(Z) \phi =
\epsilon(Z) \phi$.  Our definition of the adjoint action is chosen so
that any function $\phi$ with the property that $\phi(ZW)=\phi(WZ)$
will be invariant.  The invariant elements of $B$ form a subalgebra
denoted $B^H$.

\subsection{}

A graph consists of a set $E$ called edges, a fixed point free
involution $- :E \rightarrow E$, and a partition $V$ of $E$ into
subsets called vertices.  Let $\mathfrak{i} : E \rightarrow V$ be the
map sending $e$ to the vertex $v$ containing it.  Let $\mathfrak{t} =
\mathfrak{i} \circ -$.  We call $\mathfrak{t}(e)$ the terminal vertex
of $e$ and $\mathfrak{i}(e)$ the initial vertex.  An {\bf orientation}
is a choice $O$ of one edge from each orbit of the involution.  An
oriented graph is denoted by the data $(E,-,V,O)$.

There is a one-dimensional CW-complex associated to $(E,-,V,O)$, which
is called its {\bf geometric realization}.  The 0-cells are in one to
one correspondence with $V$, the 1-cells are in one to one
correspondence with $O$, and the characteristic maps are determined by
$\mathfrak{t}$ and $\mathfrak{i}$.

A ciliation $C$ of a graph is a linear ordering of each vertex.  The
additional data is denoted $V^c$, although we will continue to use $V$
for the partition of $E$.  A {\bf lattice} is an oriented, ciliated
graph.  The geometric realization of an oriented graph is insufficient
to support a ciliation so for a lattice we construct an oriented
surface called its {\bf envelope}.  Each vertex becomes an oriented
disk and each edge in $O$ becomes an oriented band.  The orientation
of a fattened vertex induces an orientation on its boundary, one point
of which is marked with a cilium.  Attach the band corresponding to
each $e$ to the disks (or disk) at its initial and terminal ends.  The
attaching points along the oriented boundary of each disk must be
arranged in the order given by the ciliation of the vertex.  Further
annotate the resulting surface by orienting the core of each band from
$\mathfrak{i}(e)$ to $\mathfrak{t}(e)$.

\begin{figure}
\mbox{}\hfill\epsfig{file=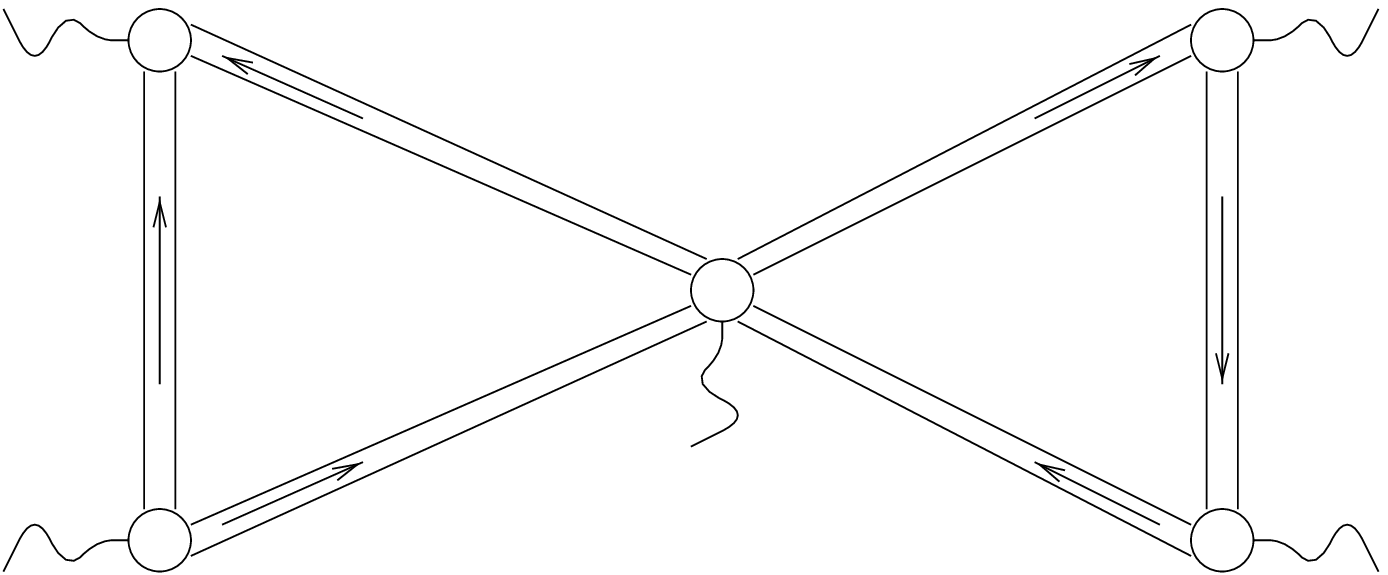,height=83pt}\hfill
\epsfig{file=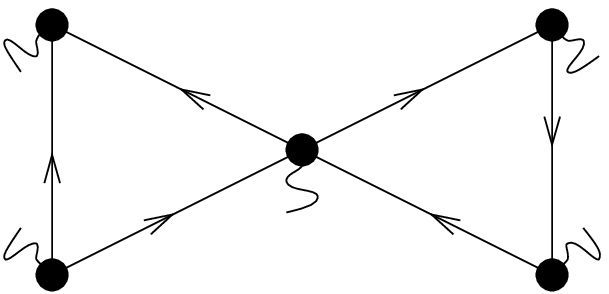}\hfill\mbox{}
\caption{Envelope of an oriented, ciliated graph.}
\label{envelope}
\end{figure}

For example, consider
\begin{align*}
E &= \{\pm e_1, \pm e_2, \pm e_3, \pm e_4, \pm e_5, \pm e_6, \}, \\ 
V^c &= \{\{-e_1,e_2\},\{-e_2,e_3\},\{-e_3,e_1,e_4,-e_6\},\{-e_5,-e_4\},
\{e_6,e_5\}\} \quad \text{and}\\
O &= \{e_1,e_2,e_3,e_4,e_5,e_6\},
\end{align*}   
Ciliation is given by the order in which the elements of each vertex
are written above.  The envelope of $(E,-,V^c,O)$ is shown in Figure
\ref{envelope}, alongsinde a streamlined schematic version.     

An envelope determines its lattice.  The edge set consists of a pair
$\pm e$ for each band, with $e$ assigned to the orientation.  Label
the initial and terminal ends of a band core by $e$ and $-e$
respectively.  Each disk forms a ciliated vertex by reading off these
labels, beginning at the cilium and traveling along the induced
orientation.

\subsection{}

For this subsection we fix $(E,-,V^c,O)$.  A gauge field theory on
this lattice is defined by the interactions of three algebraic
objects:
\begin{itemize}
\item A set of {\bf connections}, $\displaystyle{\bA =
\bigotimes_{e \in O} H}$. 
\item A {\bf gauge algebra}, $\displaystyle{\cG = \bigotimes_{v \in V}
H}$.
\item And a set of {\bf gauge fields}, $\displaystyle{C[\bA] =
\bigotimes_{e \in O} B}$.
\end{itemize}
The gauge algebra is a Hopf algebra in the natural sense of a tensor
power of Hopf algebras, whereas $\bA$ and $C[\bA]$ inherit only the
vector space structure of $H$.  However, we will shortly endow the
connections with a $\cG$-action and a comultiplication, which induce
dual structures on $C[\bA]$ via the evaluation pairing.

For each $v \in V^c$ there is a function $\mathrm{ord}_v : v
\rightarrow \mathbb{N}$ that assigns to $e \in v$ the
ordinal number corresponding to its position in the ciliation.
Connections become a left $\cG$-module under the action \[ \otimes_{v \in
V} \; y_v \cdot \otimes_{e \in O} \; x_e = \otimes_{e \in O} \; y_{
\mathfrak{i}(e) }^{ (\mathrm{ord}_{ \mathfrak{i}(e) }(e)) } \; x_e \;
S \left( y_{ \mathfrak{t}(e) }^{ (\mathrm{ord}_{ \mathfrak{t}(e) }(e))
} \right).\] This is a busy formula, even with the $|V|$-th order
summation over Sweedler notation suppressed.  For a graphical
description of the action (and the usual method of computing it) see
\cite{survey}. The gauge algebra acts adjointly on gauge fields via
$(f\cdot y)(x)=f(y\cdot x)$, so $C[\bA]$ is a right $\cG$-module.
A gauge field
$f$ is called an {\bf observable} if for all $y \in \cG$ we have $y
\cdot f = \epsilon(y) f$.  Observables form a submodule $\mathfrak{O}$
of $C[\bA]$.

\subsection{}\label{observables}

There is a construction of gauge fields which uses a direct
interpretation of the restricted dual of $H$.  Suppose that $W$ is a
finite dimensional left $H$-module.  We use $W^*$ to denote the dual
with $H$ acting on the left by $(x\cdot f)(v) = f(S(x) \cdot v)$.  The
action $(f \cdot x)(v) = f(x \cdot v)$ makes the dual into a right
$H$-module, denoted $W'$.  Compultiplication supplies a left action on
$W^* \otimes W$, namely $x \cdot (f \otimes v) = (x' \cdot f)
\otimes (x'' \cdot v)$.  Forcing the natural identification of $W^*
\otimes W$ with $\Hom(W,W)$ to be an intertwiner makes the later into a
left module as well.  In Sweedler notation the action is $(y \cdot
f)(v) = y'' \cdot f(S(y') \cdot v)$.  Now suppose that $\rho : H \to
\Hom(W,W)$ is the original representation.  I.e., $x \cdot v =
\rho(x)(v)$.  The reader may check that for any $x,y \in H$, we have
$y \cdot \rho(x) = \rho(ad(y)x)$.

A finite dimensional representation $\rho : H \to \Hom(V,V)$ is said
to be {\bf adapted to} $B$ if \[\{h \circ \rho\;|\;h \in
(\Hom(V,V))'\} \subset B.\] A {\bf coloring} of a lattice is a
labeling of each $e \in O$ by a representation adapted to $B$.  Let
$W_e$ denote the representation associated to $e$, and let $W_{-e} =
W_e^*$.  A coloring naturally associates the left $H$-module $W_v =
\bigotimes_{e \in v} W_e$ to each vertex.

Given a coloring, there is a map of right $\cG$-modules,
\[\bigotimes_{v \in V} W_v' \to \bigotimes_{e \in O} (W_{-e} \otimes
W_e)' \to \bigotimes_{e \in O} (\Hom(W_e, W_e))' \to C[\bA]\] defined
as follows.  The first stage is just reordering of the factors with
the natural distribution of primes over tensor products.  The next is
the canonical identification.  The last is composition with
$\otimes_{e \in O} \rho_e$, where the maps $\rho_e : H \to
\Hom(V_e,V_e)$ are the actual representations of the coloring. 

\begin{theorem}\label{gaugefields}
The images of these maps, taken over all colorings, add up to
$C[\bA]$.
\end{theorem}

\begin{proof}
This is evident after establishing the following claim: For each $f
\in B$ there is a finite dimensional $H$-module $W$ so that \[f \in
\{h \circ \rho\;|\;h
\in (\Hom(W,W))'\} \subset B.\]  Fix a nonzero $f \in B$.  Choose $I$
to be maximal among ideals of $H$ contained in $\ker f$.  Let $\rho :
H \to \Hom(W,W)$ be the representation induced by left multipication
of $H$ on $W = H/I$.  Define $T$ to be the linear span of the
functionals $\{y \mapsto f(xyz)\;|\; x,z \in H\}$.

Since $I$ is an ideal and it lies in the kernel of $f$, $T$ may be
thought of as a subspace of $W^*$.  Choose $\{x_1,\ldots,x_n\} \subset
H$ so that $x_1 = 1_H$ and so that, in the quotient, this is a basis
for $W$.  By evaluation, each $x_i$ is a functional on $T$.  Suppose
that, as a functional on $T$, $X = \sum a_i x_i = 0$.  For any $y,z
\in H$ we have $f(yXz)=0$.  Maximality of $I$ then implies linear
independence of $\{x_i\}$ on $T$, which means $T$ is all of $W^*$.

Choose a basis $\{f_1,\ldots,f_n\}$ for $W_f$ that is dual to
$\{x_i\}$.  let $\rho(x_j)$ be the matrix $M^j$ in the basis
$\{x_i\}$.  By duality of bases, we know that every element $y \in H$
can be written as \[y = z + \sum f_i(y) x_i\] where $z \in I$.  Hence
\[\rho(y) = \sum f_i(y) M^i.\] The $j$-$k$ entry of the matrix $\rho(y)$
is \[\sum_i f_i(y) M^i_{jk},\] which proves the assertion
\[\{h \circ \rho\;|\;h \in (\Hom(W,W))'\} \subset B.\] 

Finally, let $h_j$ be the function \[y \mapsto \sum_i f_i(y)
M^i_{j1}.\] Note that $h_j(x_i) = M^i_{j1}$.  The first column of
$M^i$ expresses $x_ix_1$ in the chosen basis.  However, $x_1 = 1_H$,
so we have $M^i_{j1} = \delta_{ji}$.  Since $h_j(I) = 0$, we have
shown that it agrees with $f_j$ on all of $H$.  This proves that $\{h
\circ \rho\;|\;h \in (\Hom(W,W))'\}$ contains a spanning set for
$T$.  In particular, it contains $f$.
\end{proof}

If $H$ is semisimple there is a way of getting an isomorphism out of
the construction above.  Restrict the colors to lie in an exhaustive
list of irreducible representations adapted to $B$, so that no
representation appears twice in the list. Once this has been done then
the map in theorem above becomes an isomorphism. This is the
definition of gauge fields used in
\cite{AGS,BR1}.

Let $\operatorname{Inv}(W)$ denote the invariant part of an
$H$-module.  Since the maps described above are all intertwiners, we
have the following characterization of observables.

\begin{corollary}
The sum over all colorings of the images of $\bigotimes_{v \in
V}\operatorname{Inv}(W_v)$ is equal to $\mathfrak{O}$.
\end{corollary} 

\section{Multitangles}

The goal of this subsection is to develop a functor between two
categories: the {\bf category of multitangles}, $\cM$, and the {\bf
category of connections}, $\cA$. The objects of $\cM$ are lattices,
and a morphism is a set of equivalence classes of tangles in
one-to-one correspondence with the vertices of its domain. An object
in $\cA$ is the set of connections on a lattice, viewed as a left
module over the gauge algebra of the lattice. The morphisms are pairs
of maps, one from the connections in its domain to the connections in
its range and the other between the gauge algebras.  The second
morphism allows us to pull back the connections on the range lattice
to a module over the domain gauge algebra.  The first map must
intertwine this action with the standard one.

A multitangle is described by a collection of diagrams in one-to-one
correspondence with the vertices of the domain lattice.  Each diagram
lies in a copy of $[0,1] \times [0,1]$ with the second factor
determining a height function on the entire collection.  Each such
collection is built by stacking the elementary diagrams described
below.  A multitangle is an equivalence class of these collections
under a relation that will be made explicit shortly.

\subsection{}  These are the elementary diagrams.
\begin{description}
\item[The identity]
The domain and range are the same lattice.  For each vertex $v$ there
is a diagram consisting of arcs with no crossings which are monotonic
with respect to the height function.  The arcs correspond, from left
to right, to the cilial ordering of $v$.  Arcs corresponding to edges
in the orientation are directed downwards, others are directed
upwards.  Figure \ref{trivial-example} shows an envelope and its
identity morphism.

The convention for ordering and directing arcs used here is standard
for all elementary diagrams.  It can also be derived from the envelope
of a lattice.  Allow the band cores to protrude into a fat vertex and
unroll it with a M\"obius transformation to the upper half plane,
cilium at infinity.

\begin{figure}
\mbox{}\hfill\epsfig{file=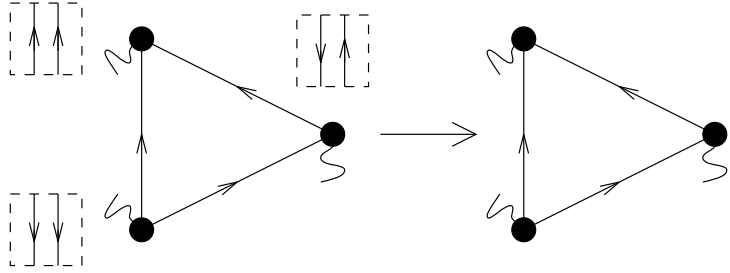}\hfill\mbox{}
\caption{Identity multitangle.}
\label{trivial-example}
\end{figure}

\item[Crossings]
The domain and range lattices differ only in the ciliation at a single
vertex, $v$, where a pair of adjacent edges have been transposed.  At
each vertex other than $v$ there is a trivial diagram as in the
identity.  The diagram at $v$ has monotonic arcs and a single crossing
between those arcs corresponding to the transposed edges.  Either
strand may pass over the other.  A geometric example is given in
Figure \ref{crossing-example}.

\begin{figure}
\mbox{}\hfill\epsfig{file=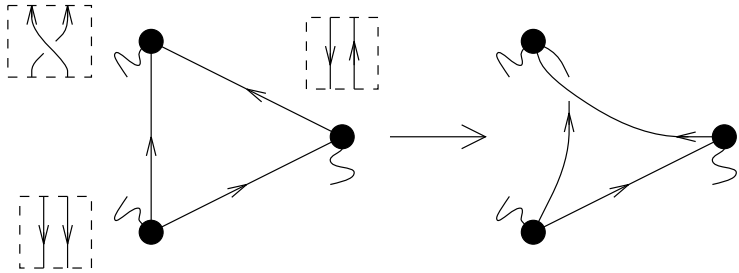}\hfill\mbox{}
\caption{Negative crossing.}
\label{crossing-example}
\end{figure}

\item[Triads]
Consider an involutary pair of edges, $e$ and $-e$, in the domain
lattice with $e \in O$.  In the range lattice, $e$ is removed from its
position in a ciliated vertex and replaced with $e',e''$ in that
order.  Similarly, $-e$ is replaced with $-e'',-e'$.  The edges $e'$
and $e''$ lie in the orientation of the new lattice, which is
otherwise identical to the domain.

The diagrams have monotonic arcs without crossings corresponding to
all the unchanged edges.  The two arcs corresponding to $e$ and $-e$
split at the same height into four arcs corresponding to $e'$, $e''$,
$-e''$ and $-e'$.  In an envelope, this is the operation of doubling an
edge (Figure \ref{triad-example}).

\begin{figure}
\mbox{}\hfill\epsfig{file=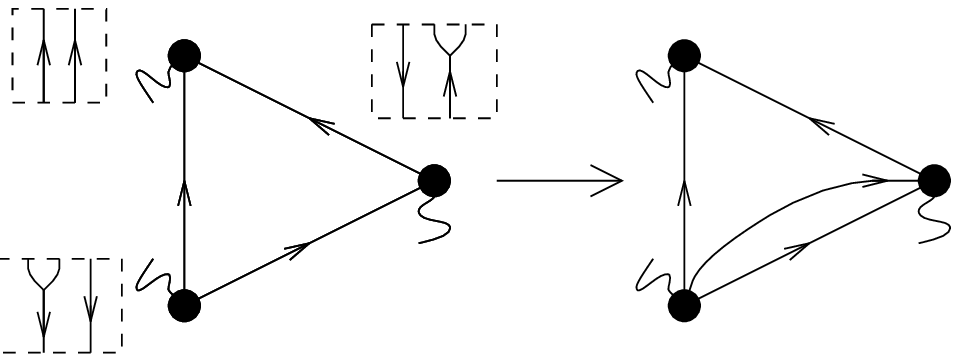}\hfill\mbox{}
\caption{Triad.}
\label{triad-example}
\end{figure}

\item[Caps] 
The edge set of the range will differ from the domain by deleting two
adjacent edges from a vertex, exactly one of which lies in the
orientation.  If the two edges are not an involutary pair, then the
two orphaned edges in the range become an involutary pair.  The
strands corresponding to the deleted edges meet at a local maximum.
Otherwise the diagrams are trivial.  The effect on envelopes is
suggested in Figure \ref{cap-example}.

\begin{figure}
\mbox{}\hfill\epsfig{file=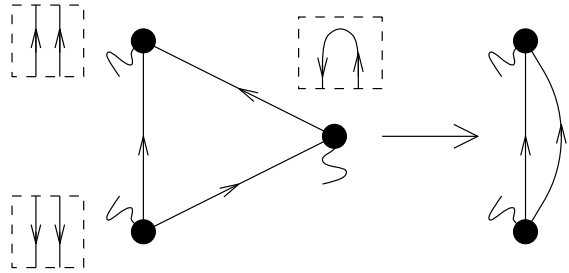}\hfill\mbox{}
\caption{Cap joining a non-involutary pair of edges.}
\label{cap-example}
\end{figure}

\item[Cups]
The range lattice differs from the domain by introducing two new
edges, $e$ and $-e$ next to each other at a single vertex.  The
strands corresponding to the new edges originate in a local minimum
and obey the usual directedness rule. Otherwise the diagrams are
trivial.  This creates a monogon at a vertex as shown in Figure
\ref{cup-example}.

\begin{figure}
\mbox{}\hfill\epsfig{file=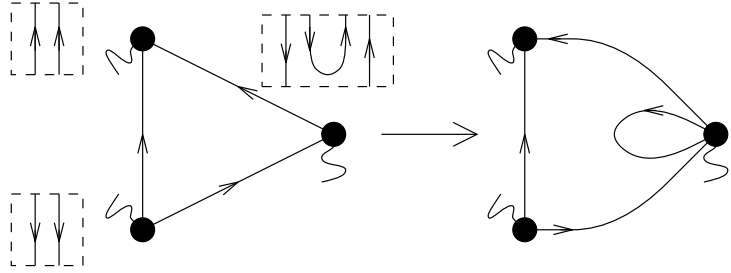}\hfill\mbox{}
\caption{Cup creating an involutary pair at a single vertex.}
\label{cup-example}
\end{figure}

\item[Stumps]
The range lattice is formed from the domain lattice by deleting an
involutary pair of edges.  The diagrams are trivial except for the two
strands corresponding to the deleted edges, which simply terminate.
Both stumps must occur at the same height.  Figure \ref{stump-example}
is an example.

\begin{figure}
\mbox{}\hfill\epsfig{file=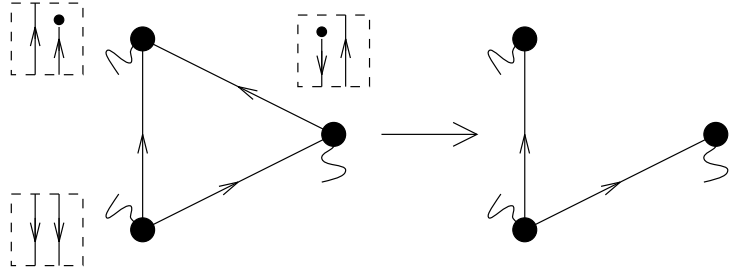}\hfill\mbox{}
\caption{Stump eliminating an involutary pair.}
\label{stump-example}
\end{figure}

\item[Switches]
The range differs from the domain by replacing some $e$ with $-e$ in
the orientation.  This is indicated in the diagrams by a hash mark on
each of the strands involved, with both marks lying at exactly the
same height  (Figure \ref{switch-example}).

\begin{figure}
\mbox{}\hfill\epsfig{file=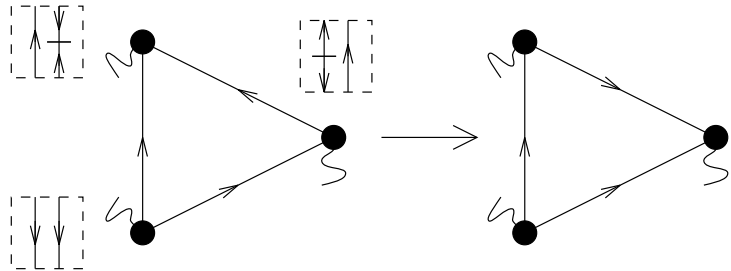}\hfill\mbox{}
\caption{Switch exchanging an involutary pair.}
\label{switch-example}
\end{figure}

\item[Cuts]
The range lattice is altered by dividing the ordered edges at some
vertex into two non-empty consecutive sets, which form new ciliated
vertices.  The diagram for that vertex is trivial, except for a
vertical mark at the top which indicates the cilium of the new vertex
(Figure \ref{cut-example}).

\begin{figure}
\mbox{}\hfill\epsfig{file=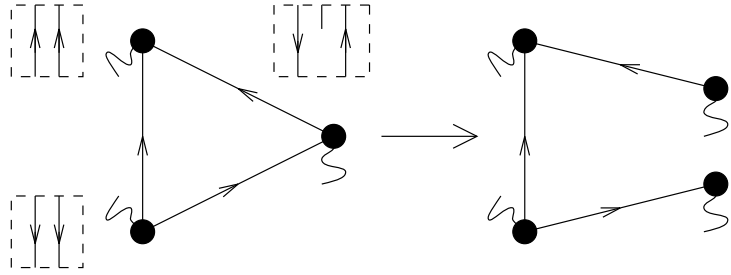}\hfill\mbox{}
\caption{Cut splitting a vertex in two.}
\label{cut-example}
\end{figure}

\end{description}

When the range lattice of one collection of diagrams matches the
domain of another one may form a new set of diagrams by stacking the 
first two.  It may be necessary to isotop the bases to get arcs to
match, and if two diagrams are stacked atop a single diagram with a
cut, the cut extends to the top of the new diagram. The height
function is then uniformly rescaled.  We define a {\bf multidiagram}
to be any such collection formed by stacking elementary diagrams.

A multidiagram is a picture of a framed embedding of a 1-dimensional
CW-complex into a collection of cubes.  The 1-cells of this complex
are called the {\bf segments} of the multidiagram.  Alternatively,
segments are the components left behind if the triads are removed and
the arcs passing under a crossing are thought of as connected.  A {\bf
coloring} of a multidiagram is an assignment of an irreducible, finite
dimensional $H$-module to each segment.  The critical points,
switches, stumps and triads of a multidiagram are collectively refered
to as {\bf events}.

\subsection{}

We say that two colored multidiagrams are equivalent if one can be
obtained from the other by a sequence of the following moves.

\begin{description}

\item[Isotopies] 
We allow ambient isotopy of the diagrams subject to the following
restricition: No two events sharing a segment may occur at the same
height.  A pair of marks indicating a stump, switch or triad must
always remain at the same height.  Cuts must remain vertical.  Triads
must maintain one segment below the horizontal and two above it.  And
the events depicted in Figure \ref{cap-moves} may not occur at the
same height if a pair of involutary edegs is represented among their
segments.

\item[Generalized Reidemeister Moves]
These are shown in Figure \ref{gen-reid-moves}.  The moves are valid
regardless of the orientations of the arcs.  We also allow the
corresponding moves with reversed crossings.

\begin{figure}
\mbox{}\hfill\epsfig{file=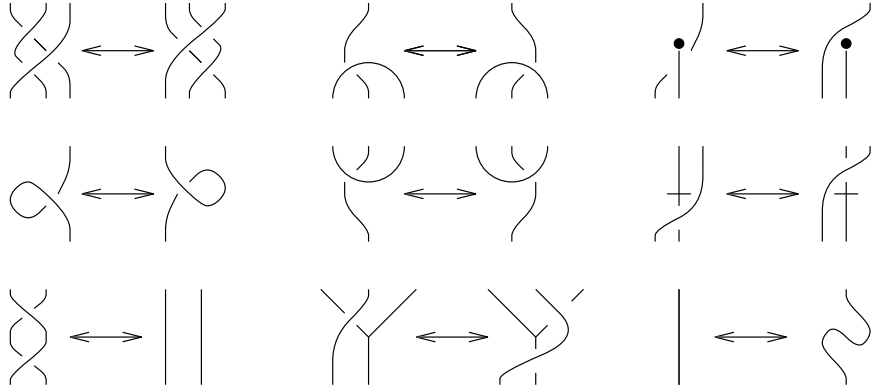}\hfill\mbox{}
\caption{Generalized Reidemeister moves.}
\label{gen-reid-moves}
\end{figure}

\item[Interacting Events]
These moves describe the interaction of events that share either a
common segment or two segments representing an involutary pair of
edges.  They are divided into triad moves (Figure \ref{triad-moves}),
cap moves (Figure \ref{cap-moves}), stump moves (Figure
\ref{stump-moves}) and switch moves (Figure \ref{switch-moves}).  In
each picture adjacent pairs of arcs, reading left to right along the
base, represent involutary pairs of edges.  Subject to that, any
orientation of arcs is allowed, as are diagrams with all crossings
reversed.  The involutary pairs are shown as adjacent merely to
conserve space; the moves are valid even for distant arcs.

Some of these moves alter the segments.  When a segment is created it
may appear with any color; a segment that splits in two takes its
color to both of the new ones; and in order for two segments to join
they must carry the same color. 

\begin{figure}
\mbox{}\hfill\epsfig{file=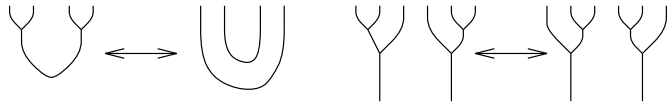}\hfill\mbox{}
\caption{Triad moves.}
\label{triad-moves}
\end{figure}

\begin{figure}
\mbox{}\hfill\epsfig{file=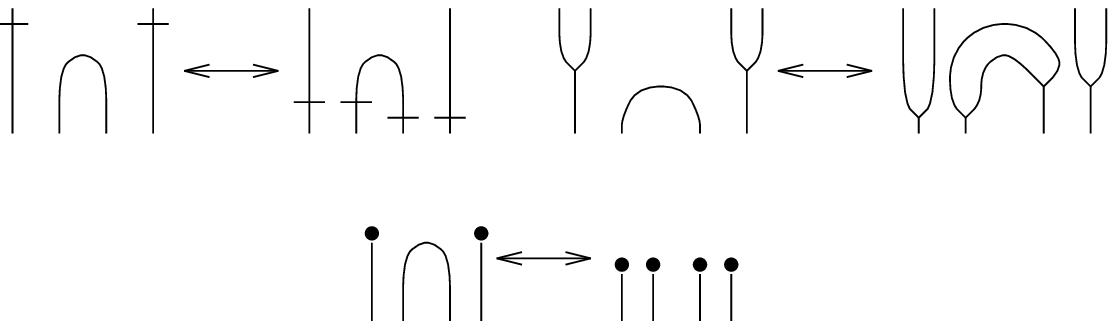}\hfill\mbox{}
\caption{Cap moves.}
\label{cap-moves}
\end{figure}

\begin{figure}
\mbox{}\hfill\epsfig{file=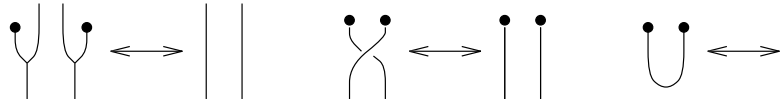}\hfill\mbox{}
\caption{Stump moves.}
\label{stump-moves}
\end{figure}

\begin{figure}
\mbox{}\hfill\epsfig{file=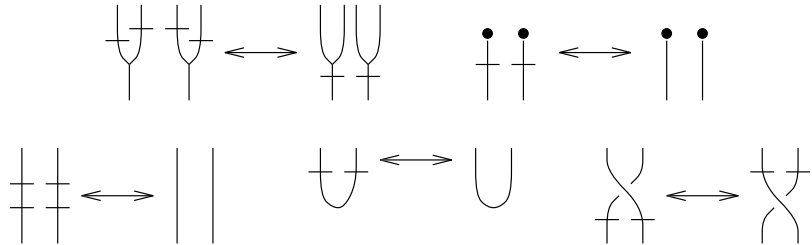}\hfill\mbox{}
\caption{Switch moves.}
\label{switch-moves}
\end{figure}

\item[Algebraic Moves] 
The two moves in Figure \ref{alg-moves} represent fundamental
identities in $H$: the definition of $S$, and $R \Delta =
\Delta^{op}R$.  As above, adjacent strands are involutary pairs and
colorings must be consistent.

\begin{figure}
\mbox{}\hfill\epsfig{file=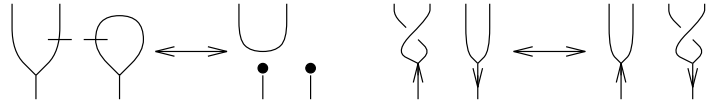}\hfill\mbox{}
\caption{Algebraic moves.}
\label{alg-moves}
\end{figure}

\end{description}

\begin{definition} A {\bf multitangle} is an equivalence class of
colored multidiagrams.
\end{definition}

There is a useful (although somewhat imprecise) topological way of
understanding the equivalence of multidiagrams.  Think of a
multidiagram as a diagram of framed tangles in cubes.  For the most
part any isotopy relative to the boundary of the cubes is an
equivalence, exceptions being the rigidities listed above.  Since
these involve involutary pairs, and since caps can alter the
involution, it is best to be careful when isotoping an event past a
cap.  Stumps are free to move almost anywhere and they are absorbed
(or created) by cups and triads.

A switch is a  pair of marks that can slide up or down together
unless obstructed by a cup, cap or triad.  Any pair of switches that
meets will cancel (and thus can be created), and a single switch can
be canceled (or created) at a cup.  A switch can move up through a
triad provided it splits in two, or two switches can combine by moving
down through a triad.

Triads, moving in pairs, can pass over each other, and a pair merging
with a cup creates a pair of cups.  Stumps can be retracted or
extended at will, and they are absorbed (or created) by cups and
triads.

In other words, as long as one avoids caps and keeps triads and stumps
upright, any continuous deformation of a multitangle is an equivalence
and the behavior of cups, stumps, triads, and switches is fairly
intuitive.  Fortunately, in practical situations caps almost always
reside above all other events in the multitangle.  If they must be
moved about, one can always rely on the list of cap moves.

In many applications the coloring of a multitangle is irrelevant.  In
those cases when it does matter, one rarely sees the moves that alter
segments.  

\subsection{}
 
We begin building a functor from $\cM$ to $\cA$ by sending a lattice
to its connections and the elementary diagrams to the morphisms
described below.

\begin{description}

\item[The identity]
This diagram induces the identity on connections and on the gauge
algebra. 

\item[Crossings]
The map on connections is an action of the $R$-matrix or its inverse.
There are 12 cases depending on the sign of the crossing, the
directions of the arcs, and the possibility that they are an
involutary pair.  These are given in Figure \ref{crossing-table},
which describes the action in the factors corresponding to the
crossing arcs.  The map extends linearly on connections, using the
identity in all other factors.  The map on the gauge algebra is the
identity.

\begin{figure}\input{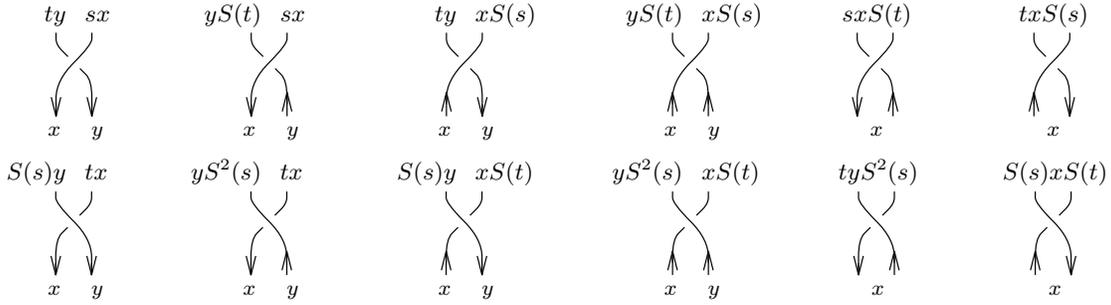}\caption{Action of $R$ and $R^{-1}$ at crossings.}\label{crossing-table}\end{figure}

\item[Triads]
Suppose that $e$ and $-e$ are replaced by $e'$, $e''$, $-e''$ and
$-e'$, with $e \in O$ originally.  The map on connections is
comultiplication in the factor corresponding to $e$ and the identity
elsewhere, with the requirement that the image of the comultiplication
take values in the tensor product of the factors corresponding to $e'$
and $e''$ in that order.  The distribution of the output, using
Sweedler notation with summation suppressed, is illustrated in Figure
\ref{triad-table}.  The diagram acts trivially on $\cG$.

\begin{figure}\input{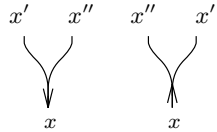}\caption{Action of a triad.}\label{triad-table}\end{figure}

\item[Caps]
There are four cases at a local maximum, depending on orientations and
on whether or not the incoming strands represent an involutary pair in
the domain lattice.  These are listed in Figure \ref{cup-cap-table},
where $\tr_V$ denotes the ordinary trace taken in the $H$-module $V$
coloring that segment.  If the lattice loses a vertex then the map on
$\cG$ is $\epsilon$ in that factor; otherwise it is the identity.

\item[Cups]
The action on $\bA$ is either by the unit of $H$ or by the unit
followed by multiplication by $k^{-1}$, depending on the orientation
of the cup.  The two cases are shown in Figure \ref{cup-cap-table}.
The map on gauge algebras is the identity.

\begin{figure}\input{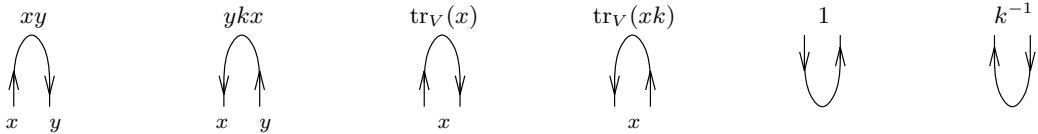}\caption{Operators corresponding to cups and caps.}\label{cup-cap-table}\end{figure}

\item[Stumps]
A stump acts as the counit in the corresponding factor of the
connections.  If the range lattice loses one or more vertices because
of this, the map on gauge algebras is counit in those factors.
Otherwise it is the identity.

\item[Switches]
For switches the map on connections is $x \mapsto S(xk)$ in the factor
corresponding to the edge.  The map on $\cG$ is trivial.

\item[Cuts] A cut has no effect on connections.  It acts trivially on $\cG$
except in the factor corresponding to the split vertex, where the map
is $\Delta : H_v \to H_{v'} \otimes H_{v''}$.  Here $v'$ denotes the
initial subset of $v$ after the cut, and $v''$ the final subset.

\end{description}

A general multidiagram is a composition of elementary ones, so it is
sent to the corresponding composition of maps.

\begin{theorem}\label{tangle-invariance}
Equivalent multidiagrams from $\Gamma$ to $\Gamma'$ induce identical
maps on connections and gauge algebras.  The map on connections
intertwines the action of $\cG_{\Gamma}$ on $\bA_{\Gamma}$ with the
one on $\bA_{\Gamma'}$ pulled back via the map on gauge algebras.
\end{theorem}

\begin{proof}. To check that equivalent multidiagrams induce the same
morphisms one must evaluate both sides of each move under every
possible arrangement of orientations and crossings.  Number the moves
in each of Figures \ref{gen-reid-moves}--\ref{alg-moves}, reading
left do right and down the page.  We will outline the identities and
manipulations in $H$ that make each move invariant on connections.
That both sides induce the same map on gauge algebras is elementary.

\begin{description}

\item[Invariance of generalized Reidemeister moves] \mbox{}
\begin{enumerate}
\item This follows from the fact that $R$ and $R^{-1}$ solve the
quantum Yang-Baxter equation.
\item Replacing any appearance of $S^2(x)$ with $kxk^{-1}$ proves
invariance in all cases.
\item This is essentially the identity  $\epsilon \otimes 1 (R) = 1 \otimes
\epsilon(R) = 1$.  For some orientations the fact that $\epsilon (S(x)) = \epsilon(x)$ is also needed.  
\item With the strand directed upwards the left hand side produces
the following morphism, where subscripts indicate successive
applications of $R$ and implied summation.
\begin{align*}
x \mapsto & x S^2(s_1) k S(t_1) \\
        = & x S(t_1 k^{-1} S(s_1))\\
        = & x S(t_1 \theta t_2 S^2(s_2) S(s_1))\\
        = & x S(\theta) S(t_1 t_2 S(s_1 S(s_2)))\\
        = & \theta x
\end{align*}
Similar computations show that, regardless of orientation, both sides
act as multiplication by $\theta$.  If the crossings are reversed the
action is by $\theta^{-1}$. 
\item This is similar to (2). 
\item Those cases that are not immediate follow from an application of
$S^2(x) = kxk^{-1}$.
\item $RR^{-1} = R^{-1}R = 1$.
\item Depending on crossings, use one of the identities $\Delta
\otimes 1 (R) = s_1 \otimes s_2 \otimes t_1t_2$ or  $1 \otimes \Delta
(R) = s_1s_2 \otimes t_2 \otimes t_1$.
\item Obvious.
\end{enumerate}

\item[Triad moves]\mbox{}
\begin{enumerate}
\item $\Delta(k^{-1}) = k^{-1} \otimes k^{-1}$.
\item $\Delta$ is coassociative. 
\end{enumerate}

\item[Cap moves]\mbox{}
\begin{enumerate}
\item That $S$ is an anti-algebra morphism suffices. 
\item $\Delta$ is an algebra morphism. 
\item $\epsilon$ is an algebra morphism and $\epsilon(k)=1$. 
\end{enumerate}

\item[Stump moves]\mbox{}
\begin{enumerate}
\item $\epsilon$ is the counit for $\Delta$.
\item $\epsilon(k^{-1}) =1$.
\item From earlier identities, $\epsilon \otimes \epsilon (R^{\pm 1}) = 1$
\end{enumerate}

\item[Switch moves]\mbox{}
\begin{enumerate}
\item $S$ is an anti-coalgebra map and $k$ is grouplike.
\item $\epsilon \circ S = \epsilon$ and $\epsilon(k)=1$.
\item This is the claim that $x \mapsto S(xk)$ is an involution.  It
follows from $S^2(x) = kxk^{-1}$.
\item $S(k)=k^{-1}$.  
\item $S \otimes S (R) = R$.
\end{enumerate}

\item[Algebraic Moves]\mbox{}
\begin{enumerate}
\item The definition of $S$: $\mu \circ S \otimes 1 \circ \Delta = \mu \circ 1
\otimes S \circ \Delta = \eta \circ \epsilon$.  
\item Constrained non-cocommutativity: $R \Delta (x) = \Delta^{op}(x)
R$. 
\end{enumerate}

\end{description}

Checking that every multitangle is a $\cG$-module intertwiner is again
a matter of checking each elementary diagram under all orientations
and crossings.  As above, we will indicate the essential identity or
manipulation on which the computation rests. 

\begin{description}
\item[Identity] This is obvious.
\item[Crossings] 
Since $\Delta$ coassociative, the identity $R \Delta (x) =
\Delta^{op}(x) R$ extends to any adjacent pair of factors in a power
of $\Delta$.  In Sweedler notation, \[y^{(1)} \otimes \cdots s y^{(i)}
\otimes t y^{(i+1)} \cdots \otimes y^{(n)} = y^{(1)} \otimes \cdots
y^{(i+1)} s \otimes y^{(i)} t \cdots \otimes y^{(n)}.\] This will
prove the intertwining of a gauge transformation by $y$ at a single
vertex.  Any other gauge transformation can be expressed as sums of
products of these.
\item[Triads] Coassociativity of $\Delta$ again.  The proof is trivial
in Sweedler notation.
\item[Caps] If the valence of the vertex is one or two, the result
follows from the defining equation for $S$.  If the valence is greater
than two, we use an extended version of the formula:
\begin{align*} 
y^{(1)} \otimes \cdots \otimes y^{(n-2)} & = y^{(1)} \otimes \cdots
S(y^{(i)}) y^{(i+1)} \cdots \otimes y^{(n)}\\ & = y^{(1)} \otimes
\cdots y^{(i)} S(y^{(i+1)}) \cdots \otimes y^{(n)}.
\end{align*}
\item[Cups] These work for pretty much the same reasons that caps do. 
\item[Stumps] Follows from the definition of $\epsilon$.
\item[Switches] $S$ is an anti-algebra map.
\item[Cuts] Coassociativity of $\Delta$.
\end{description}
\end{proof}

{\bf Remark:} We can think of \epsfig{file=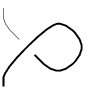,height=10pt} and
\epsfig{file=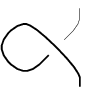,height=10pt} as single events called {\bf
positive} and {\bf negative twists} respectively.  It is worth
remembering that a positive twist acts on a connection as
multiplication by $\theta^{-1}$ in that factor.  A negative twist acts
by $\theta$.

\section{Comultiplication of Connections}

\begin{figure}
\mbox{}\hfill\epsfig{file=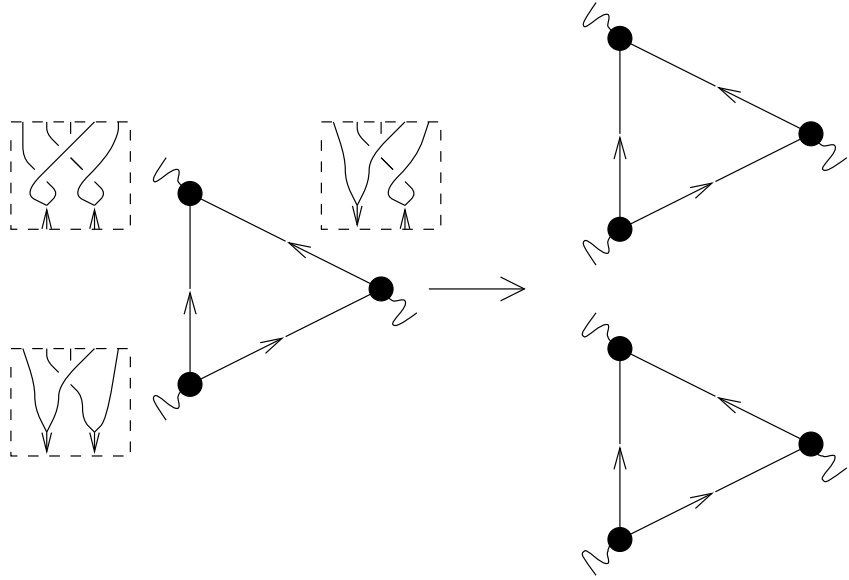}\hfill\mbox{}
\caption{Multitangle for $\nabla$.}
\label{nabla-example}
\end{figure}

Fix a lattice $\Gamma = (E,-,V^c,O)$.  We define a multitangle whose
domain is $\Gamma$ by repeating the following construction at each
vertex: Apply a triad to each arc.  Then move the strands
corresponding the the $x'$'s to the left of the the strands
corresponding to the $x''$'s so that the latter segments cross over.
Finally, cut the diagrams to separate the $x'$'s from the $x''$'s.  An
example is given in Figure \ref{nabla-example}.

The multitangle determines a range lattice denoted $\Gamma^{\otimes
2}$.  Its envelope is two disjoint copies of the envelope of $\Gamma$,
but it is important to distinguish them as the {\bf prime} and 
{\bf double prime} copies.  The induced morphism on connections is denoted 
\[\nabla : \bA_\Gamma \to \bA_{\Gamma^{\otimes 2}}.\]
The map on gauge algebras is the standard comultiplication on a
tensor power of $H$.  We make the identification $\bA_\Gamma \otimes
\bA_\Gamma = \bA_{\Gamma^{\otimes 2}}$, and define $\epsilon_\Gamma =
\otimes_{e \in O} \epsilon_e : \bA_\Gamma \to \mathfrak{b}.$

\begin{theorem}
The triple $(\bA,\nabla,\epsilon_\Gamma)$ is a coalgebra.
\end{theorem}

\begin{figure}
\mbox{}\hfill
\epsfig{file=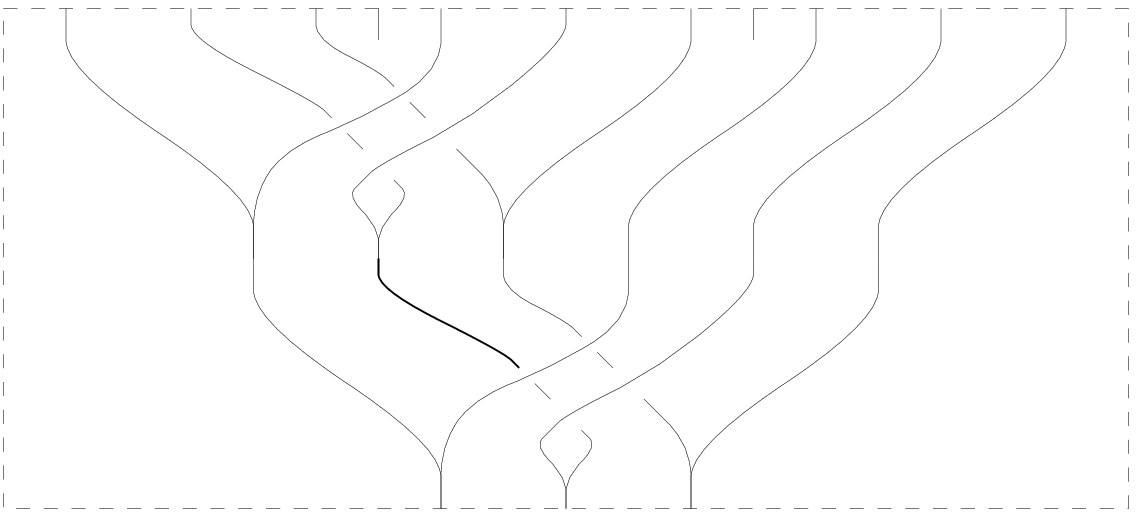,width=2in}\hfill
\epsfig{file=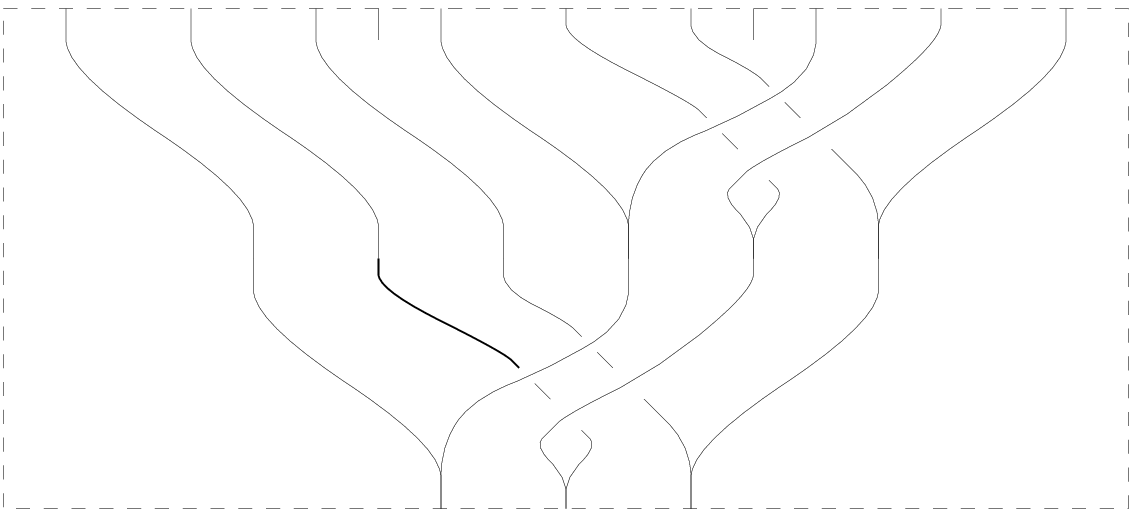,width=2in}\hfill\mbox{}
\caption{A suggestion of the coassociativity of $\nabla$}
\label{coass}
\end{figure}

\begin{proof}
Figure \ref{coass} depicts diagrams for $\nabla \otimes 1 \circ
\nabla$ and $1 \otimes \nabla \circ \nabla$ at one possible trivalent
vertex.  To see that the left side is equivalent to the right, slide
the higher triads down to the lower ones and then back up the other
segments.  This is possible because $\Delta$ is coassociative and
because the diagrams separate into three disentangled layers.  This
phenomenon holds in general, and it is possible to organize this
information into an inductive proof that $\nabla$ is coassociative.
We leave the details to the reader, with the suggestion that one use a
coupon, say
\[\raisebox{-1pt}{\epsfig{file=nabla-coupon}}\hspace{-2pt}\smlst{n}\] 
to denote the diagram for $\nabla$ at a generic $n$-valent vertex.  It
is also convenient that
\[\raisebox{-1pt}{\epsfig{file=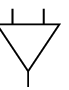}}\hspace{-2pt}\smlst{n}
\hspace{4.5pt} \raisebox{10pt}{=\;} 
\raisebox{-1pt}{\epsfig{file=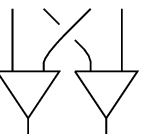}} \hspace{-26pt}
\smlst{i} \hspace{20pt} \smrst{n-i}.\] 
The fact that stumps can be retracted and absorbed into triads gives a
simple multitangle proof that $(\epsilon_\Gamma \otimes 1) \circ
\nabla = (1 \otimes \epsilon_\Gamma) \circ \nabla = 1$.  Thus
$\epsilon_\Gamma$ is a counit for $\nabla$.
\end{proof}

The adjoint of $\nabla$, restricted to observables, is denoted by
$\star$: if $f, g \in \mathfrak{O}$ and $x \in \bA$, then $(f \star
g)(x) = (f \otimes g)(\nabla(x))$.

\begin{corollary} 
$\mathfrak{O}$ is an algebra under $\star$.  
\end{corollary}

\begin{proof}
The intertwining property of morphisms induced by multitangles insures
that $\star$ takes values in $\mathfrak{O}$.  Linearity and
associativity follow from linearity and coassociativity of $\nabla$.  The
unit is the observable $\epsilon_\Gamma$.
\end{proof}

\section{Computing in $\mathfrak{O}$}  

The algebra of observables for a lattice should be independent of
orientation and should depend only on the cyclic ordering of the
ciliated vertices, not the total ordering.  Furthermore, mimicking a
classical phenomenon, the value of an observable on a connection
should be computable from a suitable connection on the complement of a
maximal tree in the graph.  In this subsection we will fix a lattice,
$(E,-,V^c,O)$, and prove that these goals can be met.  We will also
address the interaction of multitangles with the algebra and coalgebra
structures from the previous subsection.

\subsection{}  

Given $e \in O$, let $\sigma_e$ denote the map on connections induced
by the multitangle which is trivial except for a switch on the strands
$\pm e$.  In an envelope of $\Gamma$, $\sigma_e$ switches the
orientation of the core of the corresponding band.

Given $v \in V$, let
\[ \tau_v = \mbox{}^{|v|-1} \hspace{2pt}
\raisebox{-13pt}{\epsfig{file=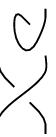}}, \quad \text{and} \quad
\tau_v^{-1} = \raisebox{-13pt} {\epsfig{ file=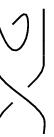}}
\hspace{2pt} \mbox{}^{|v|-1}.\]
Here an integer next to an arc indicates that many parallel
copies. The orientations are determined by the orientations of the
edges at $v$, and the rest of the multitangle is trivial.  The effect
on an envelope of $\tau_v$ is to toggle the cilium at $v$ one step
counterclockwise, while $\tau_v^{-1}$ toggles it the other way.

\begin{figure}
\centering
\mbox{}\hspace{54pt}\raisebox{-13pt}{\smst{x}}
\hspace{36pt}\smst{y_n}\raisebox{-45pt}{\smst{y_1}}\hspace{108pt}
\raisebox{-13pt}{\smst{1}}\hspace{36pt}\smst{x^{(n)}y_n}
\raisebox{-45pt}{\smst{x^{(1)}y_1}}\hspace{18pt}\mbox{}

\vspace{-45pt}

\epsfig{file=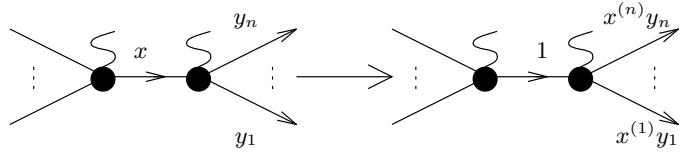}
\caption{Effect of a push on a connection}
\label{push-example}
\end{figure}

Given $e \in O$, let
\[\pi_e = \raisebox{-13pt} {\epsfig{file=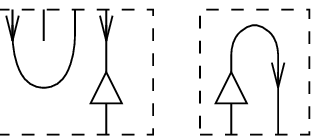}} \hspace{-54pt}
\smrst{n} \hspace{36pt} \smrst{n}\] 
where the coupon is $\Delta^{n-1}$ and the two strands entering it
represent $e$ and $-e$.  The rest of the multitangle is trivial.  The
domain and range lattices are identical.  The effect of the map on a
connection is described in Figure \ref{push-example}, which also
introduces the convention of writing a simple tensor in $\bigotimes_{e
\in O} H_e$ by labeling the corresponding cores in an envelope.  We
call this map a {\bf push}.

In order to avoid belaboring useless notation, we will assume that the
domains of successive applications of switches, toggles and pushes are
clear, provided the original domain was specified.  We will also
suppress the subscripts whenever possible.  Any sequence of toggles,
switches and pushes defines a $\cG$-module map between connection
algebras and thus an operator between observables as well.

We say that two connections are {\bf gauge equivalent} if their
difference lies in the span of $\{y \cdot x - \epsilon(y)x \;|\; x \in
\bA, y \in \cG\}$.  Observables cannot distinguish gauge equivalent
connections.  Two morphisms in $\cA$ (with the same domain and range)
are {\bf gauge equivalent} if the images of every connection are gauge
equivalent.  The adjoints of a pair of gauge equivalent operators are
identical maps on observables.

\begin{proposition}\label{cycles}
Let $f$ be any sequence of toggles and switches that begins and ends
at the same lattice.  The induced operator on connections is gauge
equivalent to the identity.
\end{proposition}

\begin{proof}
It suffices to check the compositions $\sigma^2$, $\tau^{\pm}
\tau^{\mp}$, $\sigma \tau^{\pm} \sigma \tau^{\mp}$, and $\tau_v^{\pm
n}$, where $n$ is the valence of $v$.  Clearly $\sigma$ is an
involution.  It follows easily from tangle equivalence that the next
two are also the identity map.

The multitangle for $\tau^n$ is trival away from $v$.  At that vertex
it is represented by
\[\raisebox{-10pt}{\epsfig{file=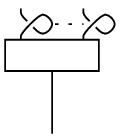}}\hspace{-11pt}\smlst{n}\]
where the coupon denotes the generator (with positive crossings) of
the center of the $n$-strand braid group.  Note that it consists
only of crossings and twists.  Because such a tangle acts as
multiplication in the factors corresponding to each segment, its
behavior can be understood by its effect on the connection 1.  Noting
that 1 is grouplike, we can evaluate this as follows:
\[\tau^n(1) = 
\raisebox{-10pt}{\epsfig{file=tauv1}}
\hspace{-11pt}\smlst{n}\hspace{11pt}= 
\raisebox{-10pt}{\epsfig{file=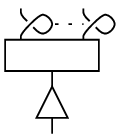}}
\hspace{-8pt}\smlst{n}\hspace{8pt}= 
\raisebox{-10pt}{\epsfig{file=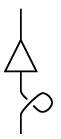}}
\raisebox{15pt}{\smlst{n}}= 
\Delta^{n-1}(\theta^{-1}).\] 
Hence, the effect on an arbitrary connection is gauge action by
$\theta^{-1}$ at $v$.  The same proof with all crossings reversed
shows that $\tau_v^{-|v|}$ is gauge action by $\theta$.
\end{proof}

\subsection{}

The standard tools for manipulating connections and observables are
toggles, switches, pushes, triads (in succession), cups, caps
(involving non-involutary pairs), stumps and cuts.  We will need an
understanding of how well they interact with the coalgebra and algebra
structures.

We already have notation for the first three maps.  Let $\Delta^n_e$
be the map induced by a succession of triads on the edges $e$ and
$-e$.  We extend the notation to include $\Delta^0$ for the identity
and $\Delta^{-1}$ for a stump.  A cup oriented from $e_1$ to $e_2$
(necessarily involutive) induces the map $\eta_{e_1}^{e_2}$.  A cap
oriented from $e_1$ to $e_2$, a non-involutive pair, induces
$\mu_{e_1}^{e_2}$.  A cut between $e_1$ and $e_2$ (in that order)
induces $\kappa_{e_1}^{e_2}$.  The notation is deliberately parallel with
the defining maps of $H$, and when context permits we will
suppress sub- and superscripts.

\begin{theorem}
The maps $\sigma$, $\Delta^n$, $\eta$ and $\mu$ are coalgebra morphisms. 
\end{theorem}

\begin{figure}
\centering
\epsfig{file=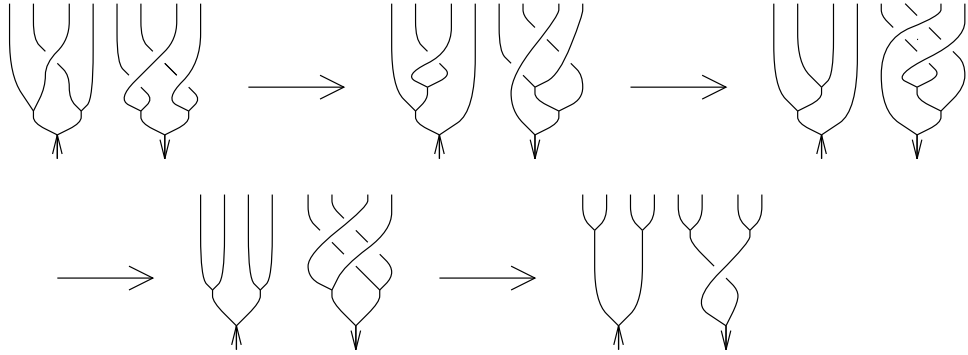}
\caption{Triad commuting with $\nabla$ at a valence two vertex.}
\label{nabla-triad}
\end{figure}

\begin{proof}
A switch slides up through the multitangle for $\nabla \circ \sigma$,
becoming a pair of switches on the appropriate edges.  Now apply the
algebraic move corresponding to $R \Delta = \Delta^{op} R$ to obtain a
tangle for $\sigma \otimes \sigma \circ \nabla$.

Consider the multitangle for $\mu \otimes \mu \circ \nabla$. At the
vertex where the cap occurs we can expand this as
\[\raisebox{-1pt}{\epsfig{file=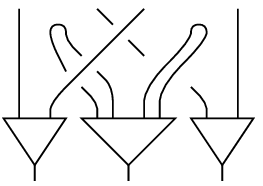}}  \hspace{-62pt}
\smrst{i} \hspace{27pt} \smrst{2} \hspace{27pt} \smrst{j}\]  
Since the caps do not involve involutary pairs, there is a cap move
that makes this into a tangle for $\nabla \circ \mu$.  The proof that
$\eta \otimes \eta \circ \nabla = \nabla \circ \eta$ is similar.

If the lattice has only two edges, the proof that $(\Delta \otimes
\Delta) \circ \nabla = \nabla \circ \Delta$ is just Figure 
\ref{nabla-triad}.  If there are more edges a similar proof works because 
of the layering phenomenon seen in the proof of Theorem 3.  Higher
powers of $\Delta$ follow immediately; $\Delta^{0}$ is trivial; and
$\Delta^{-1}$ comes from retracting stumps into triads.
\end{proof}

In many applications the output of a multitangle is needed only up to
gauge equivalence.  There is a particular occurrence which can greatly
simplify computations.  Let $\phi$ be the operator induced by a
multitangle that is trivial except for \[\raisebox{-5pt}{$\mbox{}_i$}
\hspace{-6pt} \raisebox{-13pt} {\epsfig{ file=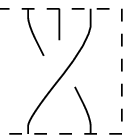}}
\hspace{-6pt} \raisebox{-5pt}{$\mbox{}_j$}\hspace{6pt}.\]  The
difference between the range and domain lattices is just that a subset
of the edges of some vertex has been split off to form a new one.
Although there is no multitangle to express it, the identity map on
connections is an operator with the same domain and range as $\phi$.

\begin{proposition}\label{reorder}
With notation as above, $\phi$ is gauge equivalent to the identity.
\end{proposition}

\begin{proof}
Given a connection $x$, we can express $\phi(x)$ as an action of
$\phi(1)$ on $x$, as in the proof of Proposition 1.  In the interest
of computing $\phi(1)$, replace the $i$ strands of the tangle with a
single strand followed by $\Delta^{i-1}$.  Tangle equivalence allows
the triads to slide over the crossing, after which we find that the
action is by $\Delta^{i-1}(s_j \cdots s_2s_1)$ in the $i$ strands and
by $t_1 \otimes t_2 \otimes \cdots \otimes t_j$ in the $j$ strands.

Now consider the effect of $\phi$ on $x_1 \otimes \cdots \otimes x_j
\otimes y_1 \otimes \cdots \otimes y_i$ in the factors corresponding 
to the non-trivial part of the tangle.  The result is
\[t_i \cdot x_1 \otimes \cdots \otimes  t_j \cdot x_j \otimes 
\Delta^{i-1}(s_j \cdots s_1) \cdot (y_1 \otimes \cdots \otimes y_i),\]
where $t_k \cdot x_k$ means $t_ky_k$ or $y_kS(t_k)$, depending on
orientation of the segment, and similarly for $\Delta^{i-1}(s_j \cdots
s_1) \cdot ( y_1 \otimes \cdots \otimes y_i)$.  This is exactly gauge
action of $1 \otimes (s_j \cdots s_1)$ on the connection 
\[t_i \cdot x_1 \otimes \cdots \otimes  t_j \cdot x_j \otimes 
y_1 \otimes \cdots \otimes y_i,\] 
which is gauge equivalent to 
\[\epsilon(s_j \cdots s_1) t_i \cdot x_1 \otimes \cdots \otimes  
t_j \cdot x_j \otimes y_1 \otimes \cdots \otimes y_i.\] 
Since $\epsilon \otimes 1 (R) = 1$, this is indistinguishable from the
behavior of the identity map.
\end{proof}

\begin{theorem}
Up to gauge equivalence, $\tau^{\pm}$, $\kappa$ and $\pi$ are
coalgebra maps.
\end{theorem}

\begin{proof}
Draw the tangle for $(\tau_v \otimes \tau_v) \circ \nabla$.  By
Proposition \ref{reorder} we can replace the cut with \[\mbox{}^{|v|}
\hspace{-12pt} \raisebox{-13pt} {\epsfig{ file=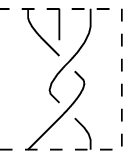}}.\] The
result is equivalent to the tangle for $\nabla \circ \tau_v$.
Commutation of $\tau^{-1}$ derives from its gauge equivalence with
some power of $\tau$.

For $(\kappa \otimes \kappa) \circ \nabla$ the tangle at the cut
vertex is
\[\epsfig{file=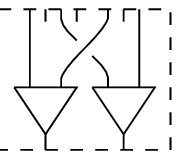}.\]  
By Proposition \ref{reorder}, this is gauge equivalent to $\nabla
\circ \kappa$.  Since $\pi$ is composed of triads, caps, a cut and a
cup, it too commutes up to gauge equivalence.
\end{proof}

\begin{corollary}
The adjoints of $\tau^{\pm}$, $\sigma$, $\pi$, $\Delta^n$, $\eta$, $\mu$
and $\kappa$ restricted to observables are algebra maps.
\end{corollary}

Since $\tau$ and $\sigma$ are invertible algebra morphisms, we can now
see that $\mathfrak{O}$ is independent of orientation, up to isomorphism, and that it depends only on the cyclic ordering of edges at vertices. 

\section{Pushes}

The goal of this section is to prove that a push is invisible to any
observable.  Let's begin with a basic fact about invariant tensors.
Suppose that $U$ and $W$ are left $H$-modules.  As in Subsection
\ref{observables}, $(U^* \otimes W)'$ is a right module.  This can be
identified with $\Hom(W,U)$ as follows: if $\phi \in U^*$, $w \in W$,
then $h \in \Hom(W,U)$ becomes functional sending $\phi \otimes w$ to
$\phi(h(w))$.  This isomorphism is an intertwiner if we use the
following action: for $z \in H$ and $h \in \Hom(W,U)$, $(h \cdot z)(w)
= S(z') \cdot h(z'' \cdot w).$

\begin{lemma}
If $h$ is invariant then, for all $z \in H$ and all $w \in W$, $h (z
\cdot w) = z \cdot h(w)$.  
\end{lemma}

\begin{proof}
Let $y =S^{-1}(z)$.  Using Sweedler notation and invariance of $h$
under the action of $y'$, we have
\begin{align*}
h(z \cdot w) &= h(S(y) \cdot w)  = h( \epsilon(y')S(y'') \cdot w) \\
             &= \epsilon(y') h (S(y'') \cdot w) = S(y') 
                \cdot h( y'' \cdot S(y''') \cdot w) \\
             &= S(y') \cdot h( \epsilon(y'') w) = S(y) \cdot h(w) \\
             &= z \cdot h(w).  
\end{align*}
\end{proof}

Now choose a coloring of the lattice with notation $\rho_e$, $W_e$,
and $W_v$ as in Subsection \ref{observables}.  Suppose that our
lattice contains the configuration in Figure \ref{push-example} and
that the vertex on the right is $v_0=\{-e_0,e_1,e_2,\ldots,e_n\}$
(with each $e_i \in O$).  Choose a gauge field of the form $f \otimes
g$, where $f \in W_{v_0}'$ and $g \in \bigotimes_{v \neq v_0} W_v'$,
and write the connection depicted in Figure \ref{push-example} as $z
\otimes x \otimes y_1 \otimes \cdots \otimes y_n$.  

The usual method of evaluation is to convert the gauge field and the
connection into tensors in $\bigotimes_{e \in E} W_e'$ and
$\bigotimes_{e \in E} W_e$, then tensor them together and contract.
(Meaning evaluate the functionals in the $W_e'$'s on the vectors in
corresponding $W_e$'s.)  Since these contractions can take place in
any order, we can focus on just those taking place between $W_e$ and
$W_e'$ for $e \in v_0$.  To see the invariance of a push, however, we
need to think of these contractions as a composition of morphisms.
   
Let $W$ denote $W_{e_1} \otimes \cdots \otimes W_{e_n}$ and let $y =
y_1 \otimes \cdots \otimes y_n$.  Apply $\rho_e$'s so that $x \in
\Hom(W_{e_0},W_{e_0})$ and $y \in \Hom(W,W)$.  We can now view $f$,
and hence $x \circ f \circ y$, as  elements of $\Hom(W,W_{e_0})$.
Using standard identifications, this becomes a tensor in $W_{e_0}
\otimes W_{e_1}^* \otimes \cdots \otimes W_{e_n}^*$.  The evaulation
of the full gauge field is now completed by contractiong $g$ with
$\alpha \otimes (x \circ f \circ y)$.

\begin{lemma} 
With notation as above, if $f$ is invariant then $x \circ f \circ y =
1 \circ f \circ (x^{(1)}y_1 \otimes \cdots \otimes x^{(n)}y_n)$.
\end{lemma}

\begin{proof}
Choose $w \in W$.  Reinterpret $(x \circ f \circ y)(w)$ as $x \in H$
acting on $f(y(w))$.  Then apply Lemma 1.
\end{proof}

In light of Corollary 1, we have established the following result:

\begin{theorem}
The adjoint of a push is the identity on observables.
\end{theorem}

\begin{corollary}\label{path-independence}
Every sequence of toggles, switches and pushes from one lattice to
another induces the same isomorphism on observables.
\end{corollary}

\begin{proof}
Since pushes induce the identity, is suffices to prove this for just
toggles and switches.  By Proposition \ref{cycles}, any two sequences
of toggles and switches are gauge equivalent.
\end{proof}

We have succeed in proving that $\mathfrak{O}$ is independent of
orientation and that it depends only on cyclic ordering at vertices.
Also, corollary \ref{path-independence} indicates that we
can evaluate observables on whatever configuration is most favorable
for the application at hand.  There is one more task.  We want to
construct (as nearly as possible) a universal description of
observables.

\section{Quantum holonomy}

Fix $\Gamma = (E,-,V^c,O)$, with envelope $F$.  We will be quantizing
the notion of holonomy along a loop in the lattice.  In the classical
world, a loop in $\Gamma$ would be just that, the image of an oriented
$S^1$ with base point.  However, in order to quantize, we need the
image to be generic and to introduce over- and under-crossings.  That
is, we need a knot diagram, not just a loop.  Also, the base point
introduces some technicalities.

We will refer to the disks in $F$ representing the elements of $V$ as
vertices and the bands as edges.  Let $\alpha$ be a proper immersion
from a disjoint, finite collection of oriented intervals into $F$, so
that endpoints map to cilia, double points lie in vertices, and the
image in each edge consists of arcs parallel to the core.  Introduce
over- and under-crossings at each interior double point.  The
resulting object is called a {\bf q-tangle}.  The special case when
the domain is a single interval is called a {\bf q-path} if the
endpoints are distinct, and a {\bf q-loop} if not.

\subsection{}

A q-tangle, $\alpha$, determines an operator, $hol_\alpha$, from $\bA$
to a tensor power of $H$ indexed by the components of $\alpha$.  It is
defined as the composition of a pair of multitangles connecting a trio
of lattices.  The first lattice is $\Gamma$; the last is determined by
the multitangles; the intermediate one comes from $\alpha$, and it is
best defined in terms of its envelope.  Its ciliated vertices are the
vertices of $F$ that $\alpha$ meets.  There is a band along each arc
of $\alpha$ in an edge of $F$, and the cores are oriented by the
orientation of $\alpha$.

The multitangle connecting the first two lattices is formed as
follows: Apply $\Delta^{m-1}$ to each pair $\pm e$, where $m$ is the
number of times $\gamma$ meets the corresponding edge of $F$.  The
range of this morphism is a lattice identical to the intermediate one,
described above, except possibly for orientation.  Continue the
multitangle by inserting switches whenever the orientations disagree.
Coloring is irrelevant.

The second multitangle is formed from the oriented tangles created by
$\alpha$ in each vertex of $F$.  Given a vertex, apply a M\"obius
transformation mapping it to the upper half plane with the cilium at
infinity.  Choose a rectangular region that contains all of the image
of $\alpha$ except for disjoint arcs running from the top edge to
infinity.  Rescale this to a unit square.  The coloring of the
resulting multitangle can be anything compatible with the first one.
The composition of the two induces the operator \[hol_\alpha : \bA \to
\bigotimes_{\text{$c$ a component of $\alpha$}} H_c.\]
One may think of envelope of the final lattice as the image of
$\alpha$, with fattened endpoints as vertices and ciliation inherited
from $F$.

\subsection{}\label{example}
 
\begin{figure}
\lmk\mbox{}
\hspace{9pt}\smst{e_5}\hspace{45pt}\raisebox{27pt}{\smst{e_4}}
\raisebox{-27pt}{\smst{e_6}}\hspace{54pt}\raisebox{27pt}{\smst{e_1}}
\raisebox{-27pt}{\smst{e_3}}\hspace{54pt}{\smst{e_2}}\hspace{4.5pt}
\mbox{}\hfill\makebox[144pt]{}\rmk
\vspace{-70pt}
\lmk\epsfig{file=bowtie.eps}\hfill\epsfig{file=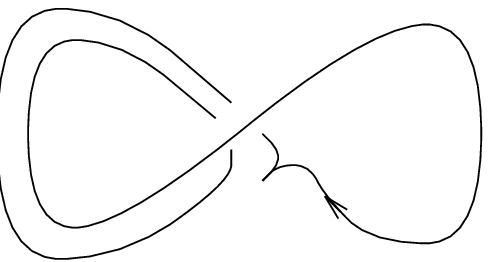}\rmk
\caption{Envelope of $\Gamma$ and a q-loop, $\alpha$.}
\label{bowtie}
\end{figure}

The notation $hol_\alpha$ is meant to be read, ``holonomy along
alpha.''  To see how it is a quantum analog of holonomy, and to
illustrate some computational devices, we will work an example where
$\Gamma$ and $\alpha$ are as in Figure \ref{bowtie}.

\begin{figure}
\lmk
\hspace{4.5pt}\smst{x_5}\hspace{49.5pt}\raisebox{27pt}{\smst{x_4}}
\raisebox{-27pt}{\smst{x_6}}\hspace{54pt}\raisebox{27pt}{\smst{x_1}}
\raisebox{-27pt}{\smst{x_3}}\hspace{49.5pt}{\smst{x_2}}\hspace{4.5pt}\rmk
\vspace{-70pt}
\lmk\epsfig{file=bowtie.eps}\rmk
\lmk\epsfig{file=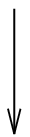}\rmk
\lmk
\smlst{S(x_5''k)}\hspace{22.5pt}\smrst{S(x_5'k)}\hspace{26pt}
\raisebox{13.5pt}{\smlst{x_4''}}\raisebox{-13.5pt}{\smlst{x_6''}}
\hspace{10pt}\raisebox{32pt}{\smrst{x_4'}}\raisebox{-32pt}{\smrst{x_6'}}
\hspace{54pt}\raisebox{27pt}{\smst{x_1}}
\raisebox{-27pt}{\smst{x_3}}\hspace{45pt}{\smst{x_2}}\hspace{4.5pt}\rmk
\vspace{-76pt}
\lmk\epsfig{file=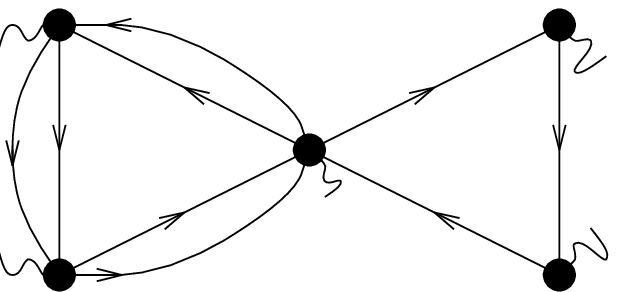}\rmk
\lmk\epsfig{file=downarrow.eps}\rmk
\lmk
\smlst{S(x_5''k)}\hspace{22.5pt}\smrst{S(x_5'k)}
\raisebox{-13.5pt}{\smrst{x_6''S^2(s_2)}}\hspace{26pt}
\raisebox{13.5pt}{\smlst{t_2x_4''}}
\hspace{10pt}\raisebox{32pt}{\smrst{t_1x_4'}}
\raisebox{-32pt}{\smrst{x_6'}}
\hspace{54pt}\raisebox{27pt}{\smst{s_1x_1}}
\raisebox{-27pt}{\smst{x_3}}\hspace{45pt}{\smst{x_2}}\hspace{4.5pt}\rmk
\vspace{-76pt}
\lmk\epsfig{file=fatbowtie.eps}\rmk
\caption{partial computation of $hol_\alpha(x)$.}
\label{hol-example}
\end{figure}

\begin{figure}
\centering
\epsfig{file=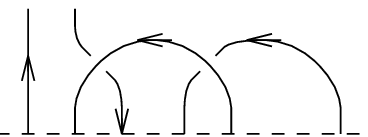}

\vspace{-6pt}

\smst{x_3}\hspace{13.5pt}\smst{x_1}\hspace{9pt}
\smst{x_4'}\hspace{18pt}\smst{x_4''}\hspace{13.5pt}\smst{x_6''}
\hspace{27pt}\smst{x_6'}
\caption{Vertex tangle.}
\label{vertex-tangle}
\end{figure}

Consider the connection $x = \bigotimes_{e_i \in O} x_i$.  Since there
is no prescribed ordering to the factors, we express $x$ by labeling
each edge of $\Gamma$ with the corresponding $x_i$.  With this
notation the evaluation of $hol_\alpha(x)$ --- except for the caps in
the multitangle --- is shown in figure \ref{hol-example}.  The top
picture is $x$.  The middle one is the result of evaluating triads and
switches.  The final stage is obtained by evaluating the crossings of
the vertex tangle shown in Figure \ref{vertex-tangle}.

We have left out evaluation of caps because there is an easier way to
do it.  Traverse the loop, concatenating the symbols accumulated on
each edge, and each time you pass through a vertex insert $k$ if the
cilium lies to your right.  This gives
\[hol_\alpha(x) = t_1 x_4' k S(x_5'' k) k x_6' k t_2 x_4'' k S(x_5'
k) k x_6'' S^2(s_2) k s_1 x_1 x_2 x_3.\] In the classical limit $k=1$
and $R =1 \otimes 1$, so $hol_\alpha$ is actually accumulating holonomy
as the loop is traversed.  Quantization occurs when self intersections
become over- or under-crossings.  The rules governing appearances of
$k$ were arrived at after lengthly experimentation. Their significance
is still unknown.

When a q-path or q-loop has no crossings it makes sense to refer to it
by listing the edges traversed.  The involution is used to indicate
the loop running against the orientation of an edge.  For example,
with $\Gamma$ and $x$ as above, we have $hol_{\{e_1,e_2,e_3\}}(x) =
x_1x_2x_3 = X$ and $hol_{\{e_4,-e_5,e_6\}}(x) = x_4S(x_5k)x_6 = Y$.
Properties of $S$ now simplify our computation to
\begin{align*}
hol_\alpha & = t_1 (x_4 S(x_5) k x_6)' k t_2 (x_4 S(x_5) k x_6)'' k
s_2 s_1 (x_1 x_2 x_3)\\ & = t_1 X' k t_2 X'' k s_2 s_1 Y
\end{align*}

In a classical setting holonomy is inverted if the direction of the
path is reversed.  In a cocommutative Hopf algebra inversion is
replaced by the antipode.  In a quantum group, however, $S$ is not an
involution.  So, we have the following quantum analog of the reversing
result for holonomy.  

\begin{proposition}\label{reversed-path}
Let $\alpha$ be a q-tangle and $\overline{\alpha}$ the same network
but with the orientation of a component $c$ reversed.  The operator
$hol_{\overline{\alpha}}$ is $hol_\alpha$ followed by a switch in the
factor indexed by $c$.
\end{proposition}

\begin{proof}
$hol_\alpha$ is defined by the composition of two multitangles, $T_1$,
consisting of stumps and triads followed by some switches, and $T_2$
formed from the vertex tangles.  If a switch is inserted between $T_1$
and $T_2$ on every edge coming from $c$, the resulting multitangle
defines $hol_{\overline{\alpha}}$.  The new switches pass upwards
through $T_2$, canceling at caps, until only one remains.  This will
lie on a pair of strands corresponding to the component $c$ in the
range of the multitangle $T_2 \circ T_1$. 
\end{proof}

\section{Wilson operators}\label{wilsonoperators}

There is an equivalence relation on q-tangles which, like the one on
multidiagrams, mimics framed tangle isotopy.  A {\bf base point free}
q-tangle is the result of smoothing the cusps at the base points of
the loops of a q-tangle.  Two base point free q-tangles are {\bf
equivalent} if they differ by a sequence of isotopies of $F$,
Reidemeister moves of type II and III, and the framing equivalence
$\raisebox{-3pt}{ \epsfig{ file=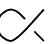}} =
\raisebox{-3pt}{ \epsfig{ file=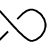}}$.  If one thinks of a base
point free q-tangle as a diagram of a framed tangle in $F \times I$,
then this equivalence is ambient isotopy fixing the endpoints of the
acr components and their framing normals.

A base point free q-tangle can be expanded into a multitangle in the
same manner as a based network.  Its evaluation, however, will depend
on the coloration.  Therefore, we define an operator for a base point
free q-tangle $L$ if and only if each closed component is colored by
an irreducible, finite dimensional $H$-module.  The segments of the
multitangle created by those components are colored accordingly; the
others receive arbitrary colors.  The induced map, denoted $W_L$,
takes values in a tensor power of $H$ indexed by the uncolored
components.

\begin{proposition}\label{loopequivalence}
If $L$ and $L'$ are equivalent base point free q-tangles with the
same coloring, then $W_L = W_{L'}$.
\end{proposition}

\begin{figure}
\centering
\epsfig{file=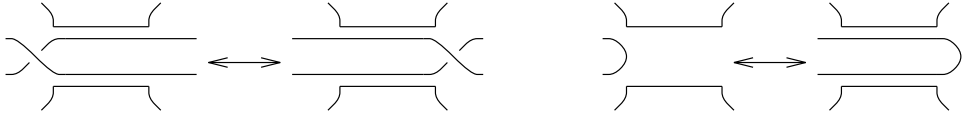}
\caption{Slides of a q-loop.}
\label{wilson-moves}
\end{figure}

\begin{proof}
Cerf theory implies $L$ and $L'$ are equivalent if, within vertices,
they differ by the moves of q-tangle equivalence, and outside vertices
they differ by the moves in figure \ref{wilson-moves}. These two moves
alter multitangles by the two algebraic moves, and the others alter
multitangles by generalized Reidemeister moves.
\end{proof}

As is usual in knot theory, we use the $L$ to denote both a base point
free q-tangle and its equivalence class.  The associated operator
$W_L$ is called a {\bf Wilson operator}.  In the special cases when
$L$ consists of closed components, a single closed component, or a
single arc component, the operators are called, respectively, a Wilson
{\bf link}, {\bf loop}, or {\bf line}.  In computing the output of a
Wilson operator, we can use the same shortcut for caps as in
\ref{example}.  The starting point on a closed component does not
matter because trace is invariant under cyclic permutation.  

Wilson operators obey a reversing result: if the reversed component is
an arc the effect is as in Proposition \ref{reversed-path}, but if it
is a closed component its statement requires more care.  Suppose that
$L$ and $\overline{L}$ differ by reversing a closed component $c$
colored by $V$.  Let $\bA_L$ denote the connections where $W_L$ and
$W_{\overline{L}}$ take values.  Define a map $f : \bA \to H_c \otimes
\bA_L$ as follows: Connect $c$ to a base point to obtain a q-tangle
$L_\bullet$.  If the cilium at which $c$ is based lies to its right,
then $f$ is $W_{L_\bullet}$ followed by right multiplication in the
factor indexed by $c$.  If it lies to the left, then $f =
W_{L_\bullet}$.  Finally, define functions $\tr_c$ and $S_c$ on $H_c
\otimes \bA_L$ as $\tr_V \otimes 1$ and $S \otimes 1$ respectively.

\begin{proposition}\label{loopreverse}
With notation as above, $W_L = \tr_c \circ f$ and $W_{\overline{L}} =
\tr_c \circ S_c \circ f$
\end{proposition}

{\bf Remark:} This is a bit easier to swallow if the $L = c$.  In that
case the Wilson loop is computed by taking the trace of something, and
the reversed loop is computed by taking $S$ of the trace of something.

\begin{proof}
Repeat the proof of Proposition \ref{reversed-path} up to the point
where all but one of the new switches have been canceled.  This time
the switch will be on two involutary strands entering a cap colored by
$V$.  The multitangles for $W_L$ and $W_{\overline{L}}$ differ by
replacing a cap with a switch followed by a reversed cap.  Since the
rest of the multitangle induces $W_{L_\bullet}$, it is easy to verify
the proposition for the two possible orientations of the cap.
\end{proof}

Since multitangles intertwine gauge transformations, Wilson links
colored by adapted representations are necessarily elements of
$\mathfrak{O}$.  Furthermore, a Wilson link is determined by a
geometric object associated to an envelope.  While toggles and
switches have no effect on the underlying q-tangle, the alterations to
the lattice will effect the evaluation of the Wilson link.
Fortunately the variance is natural.

\begin{theorem}
Suppose that $\Gamma$ and $\Gamma'$ differ by toggles and switches and
that $f$ is the induced map on connections.  If $W_L$ and $W_L'$ are
Wilson links built on the same $L$ in the two envelopes, then $W_{L'}
\circ f = W_L$.
\end{theorem}

\begin{figure}
\centering
\smst{i}\hspace{27pt}\smst{j}\raisebox{-22.5pt}{\smst{n}}\hspace{63pt}
\smst{j}\raisebox{-22.5pt}{\smst{n}}\hspace{27pt}\smst{i}
\par
\vspace{-50pt}
\par
\epsfig{file=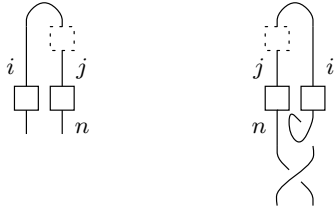}
\caption{Tangle diagrams for $W_L$ and $W_L \circ \tau_v$.}
\label{toggle-naturality}
\end{figure}

\begin{proof}
Suppose $\Gamma$ and $\Gamma'$ differ by a switch $\sigma$.  The
multitangle for $W_{L'}$ differs from that of $W_L$ by a switch on the
same edge.  Since the switches cancel, we have $W_{L'} \circ \sigma =
W_L$.  Assume now that the lattices differ by a counterclockwise
toggle, $\tau$.  That part of the multitangle for $W_L$ is represented
by the left side of Figure \ref{toggle-naturality}, where the solid
coupons are the switches and powers of $\Delta$ and the dashed coupon
contains all the crossings of the transformed vertex. The portion of
the multitangle for $W_{L'} \circ \tau$ at that vertex is now the
right side of the figure.  Since these are equivalent diagrams,
$W_{L'} \circ \tau = W_L$.  
\end{proof}

\begin{corollary}\label{natural}
Let $L$ be a colored base point free q-tangle with no arc components.
The unique isomorphism between $\mathfrak{O}_\Gamma$ and
$\mathfrak{O}_\Gamma'$ identifies the Wilson links over $L$ in each
algebra.
\end{corollary}

Observables produced by Wilson links can be graphically multiplied.
Suppose that $L$ and $L'$ are equivalence classes of base point free
q-tangles with no arcs.  Laying $L'$ over $L$ and perturbing the
result creates a new q-tangle denoted $L \star L'$.  There are several
ways to do this, but all are equivalent $q$-tangles and the result is
independent of the representatives used.

\begin{theorem}
$W_L \star W_{L'} = W_{L \star L'}$.
\end{theorem}

\begin{proof} Let $F$ denote the envelope of the lattice.  If $w_t$ is a
tangent vector along a band core, the right hand side of the edge is
determined by a normal vector $w_n$ so that $\{w_n,w_t\}$ is in the
orientation of $F$.  By equivalence of q-links, we may assume that $L$
lies to the right of $L'$ in every edge of $F$.  At each vertex, the
multitangle for $W_{L \star L'}$ begins with an arc for each edge of
$F$, oriented accordingly.  Apply a triad to each of these arcs so
that it splits into a prime branch and a double prime branch.  If $L$
meets a given edge $m$ times, apply $\Delta^{m-1}$ to the prime branch
for that edge.  If $L'$ meets it $n$ times, apply $\Delta^{n-1}$ to
the double prime branch.  Now insert the necessary switches to make
this the first half of the multitangle.

Next, focus on the second part of the multitangle at a single vertex
$v$ of $F$.  The transformed vertex fits onto the first part of the
multitangle so that the arcs of $L \cap v$ meet the prime branches and
the arcs of $L' \cap v$ meet the double primes.  Since $L$ lies under
$L'$ one can drag $L \cap v$ to the left and $L' \cap v$ to the right
until they are disjoint.  Finally, move everything between the initial
triad on each arc and $L \cap v$ to just beneath $L \cap v$, and
similarly for $L' \cap v$.  This can be done for every vertex while
preserving the simultaneous levels of triads and switches.  Inserting
a cut in between $L \cap v$ and $L' \cap v$ in each diagram, we have a
multitangle that represents $(W_L \otimes W_{L'}) \circ \nabla$. 
\end{proof}
 
Let $\mathfrak{b}\cL_F$ denote the linear space over the set of all
framed, oriented, colored links in $F \times I$ (completed if
$\mathfrak{b} = \bC[[h]]$).  This is an algebra under $\star$
which---for every $\Gamma$ whose envelope is homeomorphic to $F$---is
identified with a sub-algebra $\cW_\Gamma
\subset \mathfrak{O}_\Gamma$.  Corollary \ref{path-independence} makes
these algebras into a category.  The content of Corollary
\ref{natural} is that $\mathfrak{b}\cL_F$ behaves somewhat like a
universal object.  One of the aims of the next section is to address
how much of $\mathfrak{O}$ is generated by Wilson links.

\part{The Structure of the Algebra of Observables}

In this part of the paper we investigate the structure of the algebra
of observables for various choices of $H$ and $B$.  Suppose first that
$G$ is a connected affine algebraic group.  We will show that the ring
of $G$-characters of a free group is the algebra of observables for a
theory in which $H$ is the universal enveloping algebra of $G$ and $B$
is its coordinate ring.  Next we consider the case when $H$ is a
Drinfeld-Jimbo deformation of the universal enveloping algebra.  Here
the observables become a deformation quantization of the classical
character ring.  Following Fock and Rosly, \cite{FR}, we show that the
Poisson structure with respect to which we are quantizing is the
standard Poisson structure, as in the work of \cite{BG,Go}.  We
consider the extent to which one can generate the ring of observables
using Wilson loops.  Among the groups (and their quantum analogs) for
which this is possible is $\g$.  We conclude with a demonstration that
the observables based on $(U_h(sl_2),\ _qSL_2)$ are exactly the
Kauffman bracket skein algebra of a cylinder over the envelope of the
lattice.

\section{Classical LGFT}

\subsection{}
Let's recap a few facts from \cite{Ho}.  An affine algebraic group $G$
is a group equipped with a finitely generated algebra of
``polynomial'' functions that separate points---the coordinate ring,
denoted $B$---so that multiplication and inversion are polynomial
maps.  If $G$ is connected then its coordinate ring is naturally
identified with a stable subalgebra of its universal enveloping
algebra, which we denote by $H$.  $G$ acts on itself by conjugation,
which adjointly induces an action of $G$ on $B$.  The most important
fact for us is that the fixed part of $B$ under this action is exactly
the ring $B^H$ as defined in Subsection
\ref{hopfalgebras}.

Now suppose that $\pi$ is any finitely presented group.  A
representation of $\pi$ into $G$ is a point in $G \times \cdots \times
G$ whose coordinates are the images of the generators.  These must
satisfy polynomial identities corresponding to the relations of the
presentation, so the space of all representations is an affine
algebraic set.  It has a coordinate ring on which $G$ acts by the
adjoint of conjugation in each factor.  The fixed subring---which is
independent of presentation (up to isomorphism)---is called the {\bf
affine $G$-characters of $\pi$}.  We often shorten this to
``characters'' or ``$G$-characters'' when $\pi$ or $G$ (or both) are
understood from context.  The notation is $\cX_G(\pi)$.  If $\pi$ is
the fundamental group of a compact manifold $M$, we write $\cX_G(M)$
for $\cX_G(\pi_1(M))$.

\subsection{}
Fix a connected affine algebraic group $G$ and a lattice $\Gamma = (E,
-,V^c,O)$.  Let $H$ be the universal enveloping algebra of $G$ and $B$
its coordinate ring, thought of as lying in $H^o$.  In order to see
what the observables of this theory look like, we rebuild it using
groups instead of algebras.

Let the connections be the set of functions $A: O \rightarrow G$. The
gauge group is the set of functions $g: V \rightarrow G$, and the
gauge fields are the set $\otimes_{e \in O} B$. We can view the
connections as the Cartesian product of copies of $G$ indexed by
$O$. The gauge fields are just the coordinate ring of the Cartesian
product.  The action of the gauge group on the connections is given by
\[ g \bullet A(e)=g(\mathfrak{i}(e))A(e)g^{-1}(\mathfrak{t}(e)),\]
where $e \in O$, $g$ is an element of the gauge group, and $A$ is a
connection.  This induces a right action of the gauge group on the
gauge fields by taking adjoints.  Let $\mathfrak{O}_G$ denote the
gauge fields fixed by this action.

\begin{theorem}
Let $G$, $H$, $B$, $\Gamma = (E,-,O,V^c)$ and $\mathfrak{O}_G$ be as
above, and let $\mathfrak{O}$ denote the usual observables.
\begin{enumerate}
\item 
$\mathfrak{O} = \mathfrak{O}_G$.
\item 
$\mathfrak{O}_G$ is the $G$-characters of $\pi_1$ of the geometric
realization of $(E,-,V,O)$.
\end{enumerate}
\end{theorem}

\begin{proof}\mbox{}
\begin{enumerate}
\item 
Note that the gauge group and the gauge algebra act on the same set of
gauge fields.  A gauge field is invariant under the action of the
gauge algebra if and only if it is invariant under elements of the
form $1 \otimes \cdots \otimes y \otimes \cdots \otimes 1$, where $y$
lies in the Lie algebra of $G$.  The proof that these fields are
exactly those fixed by the gauge group action is now a simple
generalization of \cite[Corollary IV.3.2]{Ho}.
\item 
If the lattice has one vertex containing all the edges then
$\mathfrak{O}_G$ is, by definition, the $G$-characters of the
fundamental group of the geometric realization.  For arbitrary
$\Gamma$ consider another lattice $\Gamma'$ that has a single vertex,
and so that the two geometric realizations are homotopic.  Choose a
maximal tree in the geometric realization of $\Gamma$.  There is an
obvious map from the connections on $\Gamma'$ into those on $\Gamma$
that sends them to the edges of $O$ not appearing in the tree.  That
this is an isomorphism at the level of observables is an elementary
exercise.
\end{enumerate}
\end{proof}

\section{Quantized Characters}

For this subsection we assume that $(H,B)$ is defined over
$\mathbb{C}[[h]]$, and $H$ is a Drinfeld-Jimbo quantization of a
simple Lie algebra $\mathfrak{g}$ \cite{kassel}.  In this case, $H/hH$
is the universal enveloping algebra of $\mathfrak{g}$. We will denote
such $H$ by $U_h(\mathfrak{g})$.

Suppose that $B /hB$ is the coordinate ring of a connected affine
algebraic group $G$.  The equivalence of the representation theory of
$U_h(\mathfrak{g})$ and $U(\mathfrak{g})$ makes it easy to see that
$\mathfrak{O}/h\mathfrak{O}$ is the observables associated to the
theory based on $(U(\mathfrak{g}),B/hB)$. We use topological tensor
products as in \cite{kassel}, so the gauge fields are topologically
free.  Since observables form a closed subalgebra, it is also
topologically free. Therefore the observables based on
$(U_h(\mathfrak{g}),B)$ are a deformation quantization of the
$G$-characters of the fundamental group of the geometric realization
of $\Gamma$.

The only question that remains unanswered is which Poisson structure
on the ring of characters is the tangent vector to this deformation.
We can assume that $U_h(\mathfrak{g})$ has an $R$-matrix of the form
$1 \otimes 1 + h {\mathbf r} +h^2{\mathbf a}$, where ${\mathbf a}$ is
some formal power series with coefficients in $U(\mathfrak{g})\otimes
U(\mathfrak{g})$, and ${\bf r}$ is a solution of the modified
classical Yang-Baxter equation. In specific, ${\mathbf r} \in
\mathfrak{g} \otimes \mathfrak{g}$. Such solutions were classified by
Belavin and Drinfeld \cite{KS}.

The standard presentation of $U(\mathfrak{g})$ is in terms of
generators $X_i,Y_i,H_i$, where $i$ runs over some index set, and each
triple generates a copy of $U(sl_2)$. The Killing form is an element
of $\mathfrak{g}^* \otimes \mathfrak{g}^*$, but being nondegenerate,
you can contract it with respect to itself to give an element of
$\mathfrak{g} \otimes \mathfrak{g}$. This element can be expressed in
the form $a_i X_i \otimes Y_i + b_i H_i \otimes H_i + c_i Y_i \otimes
X_i$.  The element ${\mathbf r}$ can be assumed to be of the form
\[ \mathbf{r}= 2a_i X_i \otimes Y_i + b_i H_i \otimes H_i .\]
There are two actions of $\mathbf{r}$ on $B \otimes B$. You can act on
the left, by letting $\mathbf{r} \cdot (f \otimes g) (Z_1 \otimes Z_2)
= f \otimes g(\mathbf{r} \; Z_1\otimes Z_2)$, or you can act on the
right by letting $(f \otimes g) \cdot \mathbf{r} (Z_1 \otimes Z_2)=
f(Z_1 \otimes Z_2 \; \mathbf{r})$.  The Poisson bracket on $G$ is
given by \[ \{f,g\}=(f \otimes g) \cdot \mathbf{r} - \mathbf{r} \cdot
(f \otimes g).\]
 
In order to write out a formula for the Poisson structure, we need to
create a version of $\mathbf r$ that operates on the gauge fields on a
lattice. We also need to distinguish left actions and right actions:
if $Z$ acts by right multiplication in the factor corresponding to the
edge $e$, we denote it by $Z'(e)$; if it acts by multiplication on the
left, we denote it by $Z(e)$.  Since the formula we derive involves
antisymmetrization, we let $Z_1 \wedge Z_2= Z_1
\otimes Z_2 - Z_2 \otimes Z_1$.  

Now suppose that $f$ and $g$ are gauge fields in the quantum theory.
One may compute $f*g-g*f$ by writing it as $f \otimes g \circ \nabla -
g \otimes f \circ \nabla$, and expanding the $R$ matrix as a power
series in $h$. 
The linear term is

\[ \sum_{e\in O} 2a_iX_i'(e)\otimes Y_i'(e) + b_iH_i'(e)\otimes H_i'(e)
- 2a_iX_i(e)\otimes Y_i(e) - b_iH_i(e)\otimes H_i(e)+\]
\[ \sum_{v \in V}\left( -\sum_{\alpha < \beta}
 2a_iX_i'(\alpha)\wedge Y_i'(\beta) + b_iH_i'(\alpha)\wedge H_i'(\beta)\right.\]
\[-\sum_{-\alpha < -\beta}
 2a_iX_i(-\alpha)\wedge Y_i(-\beta) + b_iH_i(-\alpha)\wedge H_i(-\beta)\]
\[+\sum_{\alpha < -\beta}
 2a_iX_i'(\alpha)\wedge Y_i(-\beta) + b_iH_i'(\alpha)\wedge H_i(-\beta)\]
\[\left.+\sum_{-\alpha < \beta}
 2a_iX_i(-\alpha)\wedge Y_i'(\beta) + b_iH_i(-\alpha)\wedge H_i'(\beta)\right)\]

This is the same as the formula derived by Fock and Rosly. Hence the
Poisson structure on our algebra of classical observables is the
complex linear extension of the standard Poisson structure on the
characters of the surface with respect to the compact group.

\section{The Kauffman Bracket Skein Module}

At the end of Section \ref{wilsonoperators} we introduced the algebra
of links, $\mathfrak{b}\cL_F$, which maps naturally into the
observables for any lattice with envelope $F$.  For lattice gauge
field theory based on one of the groups $GL_n(\bC)$, $O_n(\bC)$,
$Sp_n(\bC)$, the map is onto.  This is due to a theorem of Procesi
\cite{procesi}, stating that the invariant theory for Cartesian
products of these groups is generated by traces.  It follows that
Wilson links generate observables for theories built on Drinfeld-Jimbo
deformations of these groups as well.  Furthermore, since each of
these groups has a fundamental representation in whose tensor powers
one may find all irreducible representations, a single color suffices.

The invariant theory of $SL_n(\bC)$ is a quotient of the invariant
theory for $GL_n(\bC)$, so here again Wilson links with a single color
suffice.  If $n = 2$, we can even specify the kernel of the map from
links to observables in terms of skein relations.  A space of links
divided by skein relations is a {\bf skein module}.  We will show that
the observables for a theory built on $(\U,_qSL_2)$ are the Kauffman
bracket skein module.  (Actually, there is a sign change involved in
the morphism.  This could be eliminated by redefining the skein
module, but we prefer to keep its original form.)

To obtain this result we need an explicit formula for the $R$-matrix
in the fundamental representation.  The generators of $\U$ are $X$,
$Y$ and $H$.  Let $\underline{\frac{1}{2}}$ be the vector space
spanned by $e_{-\frac{1}{2}},e_{\frac{1}{2}}$.  The standard
representation
\[ \rho : \U \rightarrow End(\V[[h]])\]
is given by 
\[\rho(X)=\begin{pmatrix} 0 & 1 \\ 0 & 0 \end{pmatrix} \qquad
\rho(Y)=\begin{pmatrix} 0 & 0 \\ 1 & 0 \end{pmatrix} \qquad
\rho(H)=\begin{pmatrix} 1 & 0 \\ 0 & -1 \end{pmatrix}. \]
One may now expand the $R$-matrix in \cite{kassel} using these
matrices and writing $e^{h/2}$ for each appearance of $q$. The trace
on $\U$ is just the trace under this representation, and takes on
values in $\bC[[h]]$.  The action of $S$ in the fundamental
representation is
\[\begin{pmatrix} a & b \\ c & d \end{pmatrix} \mapsto
  \begin{pmatrix} d & -e^{h/2}b \\ -e^{h/2}c & a \end{pmatrix},\] So
$\tr \circ S = \tr$.

In $\g$, the Cayley-Hamilton identity is:
\[ A^2 -\tr(A)A + I =0.\]

Putting the term in trace on the other side of the equation, and
multiplying by $A^{-1}$ we get a linear homogeneous equation.:

\[ A + A^{-1} = \tr(A)I.\]

Now multiply by $B$ to get back  a bilinear equation:

\[ AB +  A^{-1}B =\tr(A)B.\]

Finally take the trace:

\[ \tr(AB) + \tr( A^{-1}B) = \tr(A)\tr(B). \]

This formula is the Cayley-Hamilton identity as a trace identity for
$\g$. It is fundamental for $\g$ in the sense that
every other identity between traces follows from evaluation of this
one.

The same identity persists in $U(sl_2)$, except that you need to use
the antipode instead of the inverse:
\[ \tr(ZW) + \tr(S(Z)W) = \tr(Z)\tr(W).\]
Finally the fundamental trace identity for $\U$ is:
\[ t\;\tr(ZW) + t^{-1}\;\tr(S(Z)W) = \sum_i \tr(s_iZ)\tr(t_iW),\]
where $\sum_i s_i \otimes t_i$ is the $R$-matrix for $\U$ and $t = e^{h/4}$. 

Let $\cL_M$ be the set of framed links (including $\emptyset$) in a
$3$-manifold $M$.  Let $\bC\cL_M[[h]]$ denote formal power series in
$h$ with coefficients in the vector space over $\cL_M$.  We define
$S(M)$ to be the topological submodule of $\bC\cL_M[[h]]$ generated by
all expressions of the form
\begin{enumerate}
\item $\displaystyle{\lcr+t\zer+t^{-1}\ift}$, \quad and
\vspace{5pt}
\item $\bigcirc+t^2+t^{-2}$.
\end{enumerate}
These formulas indicate relations that hold among links which can be
isotoped in $M$ so that they are identical except in the neighborhood
shown.  The Kauffman bracket skein module is the quotient \[K(M) =
\bC\cL_M[[h]]/S(M).\] The Kauffman bracket skein module of $F \times
I$ is an algebra with multiplication as in the space of links in
Section \ref{wilsonoperators}.  We denote it $K(F)$. From \cite{quant}
we know that $K(F)$ is a deformation quantization of the
$\g$-characters of the fundamental group of $F$.

\begin{theorem}
Let $F$ be a compact, connected surface with boundary.  If the surface
underlying the envelope of $\Gamma$ is $F$ then the observables of a
lattice gauge field theory based on $(\U,_qSL_2)$ are canonically
isomorphic to $K(F)$.
\end{theorem}

\begin{proof}
We define a map $\zeta : K(F) \rightarrow \mathfrak{O}$ as follows.
Let $L$ be a framed link in $F \times I$.  Represent it as a diagram
in the envelope of $\Gamma$ with the blackboard framing.  Orient the
components arbitrarily and color them  with the fundamental
representation.  Now perturb it so that it is a base point free
q-tangle, also denoted $L$.  Finally, map it to the observable
$(-1)^{|L|}W_L$.  This is well defined at the level of $\cL_{F \times
I}$ by Propositions \ref{loopequivalence} and \ref{loopreverse} and
the fact that $\tr \circ S = \tr$.  That $\zeta$ sends elements of
$S(F \times I)$ to zero follows from the quantum Cayley-Hamilton
identity.

By Theorem 8,
 $\zeta$ is an algebra map.  To see that it is an
isomorphism, consider the map induced between $K(F)/hK(F)$ and
$\mathfrak{O}/h\mathfrak{O}$.  This is known to be an isomorphism
(\cite{isomorphism}), so the $\zeta$ must be one  as well.
\end{proof}

\end{document}